\documentclass[a4paper]{article}
\usepackage{amsfonts}
\usepackage{amssymb}
\usepackage{amsthm}
\usepackage{amsmath}
\usepackage{latexsym}
\usepackage[T1]{fontenc}
\usepackage[latin1]{inputenc}
\usepackage{color}

\def \eop {\hbox{}\nobreak\hfill
\vrule width 2mm height 2mm depth 0mm
\par \goodbreak \smallskip}

\def \eop {\hbox{}\nobreak\hfill \vrule width 2.0mm height 1.8mm depth 0mm
\par \goodbreak \smallskip}
\numberwithin{equation}{section}

\newtheorem{definition}{Definition}[section]
\newtheorem{theorem}{Theorem}[section]

\newtheorem{lemma}{Lemma}[section]
\newtheorem{remark}{Remark}[section]

\def \eop {\hbox{}\nobreak\hfill
\vrule width 2mm height 2mm depth 0mm
\par \goodbreak \smallskip}

\input{tcilatex}
\begin{document}

\author{Badreddine Mansouri, Mostapha abd elouahab Saouli
\\
{University Mohamed Khider, PO BOX 145, 07000 Biskra, Algeria}}
\title{Reflected Discontinuous {Backward Doubly Stochastic Differential
Equation With Poisson Jumps.} }
\date{}
\maketitle


\renewcommand{\thefootnote}{
\fnsymbol{footnote}} \footnotetext{{\scriptsize E-mail: mansouri.badreddine@gmail.com
(Badreddine\ Mansouri), saoulimoustapha@yahoo.fr (Mostapha abd el ouahab Saouli)}}


\noindent \textbf{Abstract.} In this paper{\ }we prove the existence of a
solution for reflected backward doubly stochastic differential equations
with poisson jumps (RBDSDEPs) with one continuous barrier where the
generator is continuous and also we study the RBDSDEPs with a linear growth
condition and left continuity in $y$ on the generator. By a comparison
theorem established here for this type of equation we provide a minimal or a
maximal solution to RBDSDEPs.

\bigskip


\section{Introduction.}

\qquad A new kind of backward stochastic differential equations was
introduced by Pardoux and Peng [11] in 1994, which is a class of
backward doubly stochastic
differential equations (BDSDEs for short)
\begin{equation*}
Y_{t}=\xi +\int_{t}^{T}f(s,Y_{s},Z_{s})ds+\int_{t}^{T}g(s,Y_{s},Z_{s})d%
\overleftarrow{B}_{s}-\int_{t}^{T}Z_{s}dW_{s},\ 0\leq t\leq T,
\end{equation*}%
where $\xi $ is a random variable termed the terminal condition, $f:\Omega
\times \lbrack 0,T]\times
\mathbb{R}
\times
\mathbb{R}
^{d}\rightarrow
\mathbb{R}
$, $g:\Omega \times \lbrack 0,T]\times
\mathbb{R}
\times
\mathbb{R}
^{d}\rightarrow
\mathbb{R}
$ are two jointly measurable processes, $W$ and $B$ are two mutually
independent standard Brownian motion, with values, respectively in $%
\mathbb{R}
^{d}$ and $%
\mathbb{R}
$. Several authors interested in weakening this assumption see Bahlali et al [3], Boufoussi et al. [5], Lin. Q
[8] and [9], N'zi el al. [10], Shi et al. [12], Wu et al.
[14], Zhu et al. [16]. A class of backward doubly
stochastic differential equations with jumps was study by Sun el al. [13],  Zhu et al. [15] They have proved the existence and
uniqueness of solutions for this type of BDSDEs under uniformly Lipschitz
conditions.

In addition, Bahlali et al [2] prove the existence and uniqueness of
solutions to reflected backward doubly stochastic differential equations
(RBDSDEs) with one continuous barrier and uniformly Lipschitz coefficients.
The existence of a maximal and a minimal solution for RBDSDEs with
continuous generator is also established.

In this paper, we study the now well-know reflected backward doubly
stochastic differential equations with jumps (RBDSDEPs for short):%
\begin{equation}
Y_{t}=\xi +\int_{t}^{T}f(s,\Lambda _{s})ds+\int_{t}^{T}g(s,\Lambda _{s})d%
\overleftarrow{B}_{s}+\int_{t}^{T}dK_{s}-\int_{t}^{T}Z_{s}dW_{s}-%
\int_{t}^{T}\int_{E}U_{s}\left( e\right) \tilde{\mu}\left( ds,de\right) ,\
0\leq t\leq T,
\end{equation}%
where $\Lambda _{s}=\left( Y_{s},Z_{s},U_{s}\right) $.

Motivated by the above results and by the result introduced by Fan. X, Ren.
Y [6] and Zhu, Q., Shi, Y [15, 16], we establish firstly the existence of
the solution of the reflected BDSDE with Poisson jumps (RBDSDEP in short)
under the continuous coefficient, also we prove the existence solution of a
RBDSDEP where the coefficient $f$ satisfy a linear growth and left
continuity in $y$ conditions on the generator of this type of equation.

The organization of the paper is as follows. In Section 2, we give some
preliminaires and we consider the spaces of processus also we define the Itô%
's formula. In Section 3, we proof a comparison theorem, section 4 under a
continuous conditions on $f$ we obtain the existence of a minimal solution
of RBDSDEP, and finally in section 5, we study RBDSDEP where the generator $f$ satisfied a left
continuity in $y$ and linear growth conditions.

\section{Notation, assumptation and definition.}

\qquad Let $\left( \Omega ,\mathcal{F},P\right) $ be a complete
probability space. For $T>0$, We suppose that $\left(
\mathcal{F}_{t}\right) _{t\geq 0}$ is generated by the following
three mutually independent processes:

(i) Let $\left\{ W_{t},0\leq t\leq T\right\} $ and $\left\{
B_{t},0\leq
t\leq T\right\} $ be two standard Brownian motion defined on $\left( \Omega ,%
\mathcal{F},P\right) $ with values in $\mathbb{R}^{d}$ and
$\mathbb{R}$, respectively.

(ii) Let random Poisson measure $\mu $ on $E\times
\mathbb{R}
_{+}$ with compensator $\nu \left( dt,de\right) =\lambda \left( de\right) dt$%
, where the space $E=%
\mathbb{R}
-\left\{ 0\right\} $ is equipped with its Borel field $\mathcal{E}$
such that $\left\{ \tilde{\mu}\left( \left[ 0,t\right] \times
A\right) =\left( \mu -\nu \right) \left( \left[ 0,t\right] \times
A\right) \right\} $ is a martingale for any $A\in \mathcal{E}$
satisfying $\lambda \left( A\right) <\infty $. $\lambda $ is a
$\sigma $ finite measure on $\mathcal{E}$\ and satisfies
$\int_{E}\left( 1\wedge \left\vert e\right\vert ^{2}\right) \lambda
\left( de\right) <\infty .$

Let \ $\mathcal{F}_{t}^{W}:=\sigma (W_{s};0\leq s\leq t)$,
$\mathcal{F} _{t}^{\mu }:=\sigma (\mu _{s};0\leq s\leq t)$\ and \
$\mathcal{F} _{t,T}^{B}:=\sigma (B_{s}-B_{t};t\leq s\leq T),$
completed with $P$-null
sets. We put, \ $\mathcal{F}_{t}:=\mathcal{F}_{t}^{W}\vee \mathcal{F}%
_{t,T}^{B}\vee \mathcal{F}_{t}^{\mu }$. It should be noted that $\left(
\mathcal{F}_{t}\right) $ is not an increasing family of sub $\sigma -$%
fields, and hence it is not a filtration.

\noindent For $d\in
\mathbb{N}
^{\ast },$ $\left\vert \cdot \right\vert $ stands for the Euclidian norm in $%
\mathbb{R}
^{d}\times \left[ 0,T\right] .$

We consider the following spaces of processus:

\begin{itemize}
\item \vskip0.15cm We denote by $\mathcal{S}^{2}\left( 0,T,%
\mathbb{R}
^{d}\right) $, the set of continuous $\mathcal{F}_{t}-$measurable processes $%
\left\{ \varphi _{t};t\in \left[ 0,T\right] \right\} $, which satisfy $%
\mathbb{E}(\sup_{0\leq t\leq T}\left\vert \varphi _{t}\right\vert
^{2})<\infty $.

\item \vskip0.15cm Let $\mathcal{M}^{2}\left( 0,T,%
\mathbb{R}
^{d}\right) $ denote the set of \ $d-$ dimensional, $\mathcal{F}_{t}-$%
measurable processes $\left\{ \varphi _{t};t\in \left[ 0,T\right] \right\} $%
, such that \ $\mathbb{E}\int_{0}^{T}\left\vert \varphi _{t}\right\vert
^{2}dt<\infty $.

\item \vskip0.15cm $\mathcal{A}^{2}$ set of continuous, increasing, $%
\mathcal{F}_{t}-$measurable process $K:\left[ 0,T\right] \times \Omega
\rightarrow \lbrack 0,+\infty ($ with $K_{0}=0,$ $\mathbb{E}\left(
K_{T}\right) ^{2}<+\infty .$

\item \vskip0.15cm $\mathbb{L}^{2}$ set of $\mathcal{F}_{T}$- measurable
random variables $\xi :\Omega \rightarrow
\mathbb{R}
$ with $\mathbb{E}\left\vert \xi \right\vert ^{2}<+\infty .$

\item \vskip0.15cm We denote by $\mathcal{L}^{2}\left( 0,T,\tilde{\mu},%
\mathbb{R}
^{d}\right) $, the space of mappings $U:\Omega \times \left[ 0,T\right]
\times E\rightarrow
\mathbb{R}
^{d}$ which are $\mathcal{P}\otimes \mathcal{E}$ measurable such that%
\begin{equation*}
\left\vert \left\vert U_{t}\right\vert \right\vert _{\mathcal{L}^{2}\left(
0,T,\tilde{\mu},%
\mathbb{R}
^{d}\right) }^{2}=\mathbb{E}\int_{0}^{T}\left\vert \left\vert
U_{t}\right\vert \right\vert _{L^{2}\left( E,\mathcal{E},\lambda ,%
\mathbb{R}
^{d}\right) }^{2}dt<\infty ,
\end{equation*}

where $\mathcal{P}\otimes \mathcal{E}$ denoted the $\sigma -$algebra of $%
\mathcal{F}_{t}-$predectable sets of $\Omega \times \left[ 0,T\right] $ and%
\begin{equation*}
\left\vert \left\vert U_{t}\right\vert \right\vert _{L^{2}\left( E,\mathcal{E%
},\lambda ,%
\mathbb{R}
^{d}\right) }^{2}=\int_{E}\left\vert U_{t}\left( e\right) \right\vert
^{2}\lambda \left( de\right) .
\end{equation*}

\item Notice also the space $\mathcal{D}^{2}\left(
\mathbb{R}
\right) =\mathcal{S}^{2}\left( 0,T,\mathbb{R}
\right) \times \mathcal{M}^{2}\left( 0,T,
\mathbb{R}
^{d}\right) \times \mathcal{L}^{2}\left( 0,T,\tilde{\mu},
\mathbb{R}
\right) \times \mathcal{A}^{2}$ endowed with the norm
\begin{equation*}
\left\vert \left\vert \left( Y,Z,U,K\right) \right\vert \right\vert _{\mathcal{D}
^{2}}=\left\vert \left\vert Y\right\vert \right\vert _{\mathcal{S}
^{2}}+\left\vert \left\vert Z\right\vert \right\vert _{\mathcal{M}
^{2}}+\left\vert \left\vert U\right\vert \right\vert _{\mathcal{L}
^{2}}+\left\vert \left\vert K\right\vert \right\vert _{\mathcal{A}^{2}}.
\end{equation*}
\end{itemize}
is a Banach space.
\begin{definition}
A solution of a reflected BDSDEPs is a quadruple of processes
$\left( Y,Z,K,U\right) $ wich satisfies
\begin{equation*}
\left\{
\begin{array}{l}
i)Y\in \mathcal{S}^{2}\left( 0,T,%
\mathbb{R}
\right) ,\text{ }Z\in \mathcal{M}^{2}\left( 0,T,%
\mathbb{R}
^{d}\right) \text{ },K\in \mathcal{A}^{2}\text{ },U\in \mathcal{L}^{2}\left(
0,T,\tilde{\mu},%
\mathbb{R}
\right) , \\
\\
ii)\text{ }Y_{t}=\xi
+\int_{t}^{T}f(s,Y_{s},Z_{s},U_{s})ds+\int_{t}^{T}g(s,Y_{s},Z_{s},U_{s})d%
\overleftarrow{B}_{s} \\
\\
+\int_{t}^{T}dK_{s}-\int_{t}^{T}Z_{s}dW_{s}-\int_{t}^{T}\int_{E}U_{s}\left(
e\right) \tilde{\mu}\left( ds,de\right) ,\ 0\leq t\leq T, \\
\\
iii)\text{ }S_{t}\leq Y_{t},\text{\quad }0\leq t\leq T\text{\quad and\quad }%
\int_{0}^{T}\left( Y_{t}-S_{t}\right) dK_{t}=0.%
\end{array}%
\right.
\end{equation*}
\end{definition}

We give the following assumptions \textbf{(H)} on the data $\left( \xi
,f,g,S\right) $:

\vskip0.1cm\noindent \textbf{(H.1) }$f:[0,T]\times \Omega \times
\mathbb{R}
\times
\mathbb{R}
^{d}\times \mathcal{L}^{2}\left( 0,T,\tilde{\mu},%
\mathbb{R}
\right) \rightarrow
\mathbb{R}
;$ $g:[0,T]\times \Omega \times
\mathbb{R}
\times
\mathbb{R}
^{d}\times \mathcal{L}^{2}\left( 0,T,\tilde{\mu},%
\mathbb{R}
\right) \rightarrow
\mathbb{R}
$ be jointly measurable such that for any $\left( y,z,u\right) \in
\mathbb{R}
\times
\mathbb{R}
^{d}\times \mathcal{L}^{2}\left( 0,T,\tilde{\mu},%
\mathbb{R}
\right) $%
\begin{eqnarray*}
f(\cdot ,\omega ,y,z,u) &\in &\mathcal{M}^{2}\left( 0,T,%
\mathbb{R}
\right) , \\
g(\cdot ,\omega ,y,z,u) &\in &\mathcal{M}^{2}\left( 0,T,%
\mathbb{R}
\right) .
\end{eqnarray*}%
\vskip0.1cm\noindent \textbf{(H.2) }There exist constant $C> 0$
and a constant $0<\alpha <1$ such that for every $\left( \omega
,t\right) \in \Omega \times \lbrack 0,T]$ and $\left( y,y^{^{\prime
}}\right) \in
\mathbb{R}
^{2},$ $\left( z,z^{^{\prime }}\right) \in \left(
\mathbb{R}
^{d}\right) ^{2},$ $\left( u,u^{^{\prime }}\right) \in \left( \mathcal{L}%
^{2}\left( 0,T,\tilde{\mu},%
\mathbb{R}
\right) \right) ^{2}$%
\begin{equation*}
\left\{
\begin{array}{l}
\left( i\right) \text{ }\left\vert f(t,\omega ,y,z,u)-f(t,\omega
,y^{^{\prime }},z^{^{\prime }},u^{^{\prime }})\right\vert ^{2}\leq
C\left[ \left\vert y-y^{^{\prime }}\right\vert ^{2}+\left\vert
z-z^{^{\prime
}}\right\vert ^{2}+\left\vert u-u^{^{\prime }}\right\vert ^{2}\right] , \\
\left( ii\right) \text{ }\left\vert g(t,\omega ,y,z,u)-g(t,\omega
,y^{^{\prime }},z^{^{\prime }},u^{^{\prime }})\right\vert ^{2}\leq
C\left\vert y-y^{^{\prime }}\right\vert ^{2}+\alpha \left(
\left\vert z-z^{^{\prime }}\right\vert ^{2}+\left\vert u-u^{^{\prime
}}\right\vert
^{2}\right) .%
\end{array}%
\right.
\end{equation*}%
\vskip0.1cm\noindent \textbf{(H.3) }The terminal value\textbf{\ }$\xi $ be a
given random variable in $L^{2}$.\newline
\vskip0.1cm\noindent \textbf{(H.4) }$\left( S_{t}\right) _{t\geq 0},$ is a
continuous progressively measurable real valued process satisfying%
\begin{equation*}
\begin{tabular}{lll}
$\mathbb{E}\left( \sup_{0\leq t\leq T}\left( S_{t}^{+}\right) ^{2}\right)
<+\infty ,$ & where & $S_{t}^{+}:=\max \left( S_{t},0\right) .$%
\end{tabular}%
\end{equation*}%
\vskip0.1cm\noindent \textbf{(H.5) }$S_{T}\leq \xi $, $\mathbb{P}$-almost
surely.

\begin{theorem}\textbf{[6] }
\label{th0} Assume that $\left( H.1\right) -\left( H.5\right) $
holds. Then Eq $\left( 1.1\right) $ admits a unique solution $\left(
Y,Z,U,K\right) \in \mathcal{D}^{2}\left(
\mathbb{R}
\right) $.
\end{theorem}

The result depends on the following extension of the well-krown Itô's
formula. Its proof follows the same way as lemma 1.3 of [11]

\begin{lemma}
\label{Lemma 1} Let $\alpha \in \mathcal{S}^{2}\left(0,T,
\mathbb{R}
^{k}\right) ,$ $\left( \beta ,\gamma \right) \in \left( \mathcal{M}%
^{2}\left(
\mathbb{R}
^{k}\right) \right) ^{2},$ $\eta \in \mathcal{M}^{2}\left(
\mathbb{R}
^{k\times d}\right) $ and $\sigma \in \mathcal{L}^{2}\left( 0,T,\tilde{\mu},%
\mathbb{R}
^{k}\right) $ such that:%
\begin{equation*}
\alpha _{t}=\alpha _{0}+\int_{0}^{t}\beta _{s}ds+\int_{0}^{t}\gamma
_{s}dB_{s}+\int_{0}^{t}\eta
_{s}dW_{s}+\int_{0}^{t}dK_{s}+\int_{0}^{t}\int_{E}\sigma _{s}\left( e\right)
\tilde{\mu}\left( ds,de\right) ,
\end{equation*}%
then $\left( i\right) $%
\begin{eqnarray*}
\left\vert \alpha _{t}\right\vert ^{2} &=&\left\vert \alpha _{0}\right\vert
^{2}+2\int_{0}^{t}\langle \alpha _{s},\beta _{s}\rangle
ds+2\int_{0}^{t}\langle \alpha _{s},\gamma _{s}\rangle
dB_{s}+2\int_{0}^{t}\langle \alpha _{s},\eta _{s}\rangle
dW_{s}+2\int_{0}^{t}\langle \alpha _{s},dK_{s}\rangle  \\
&&+2\int_{0}^{t}\int_{E}\langle \alpha _{s-},\sigma \left( e\right) \tilde{\mu%
}\left( ds,de\right) \rangle -\int_{0}^{t}\left\vert \gamma _{s}\right\vert
^{2}ds+\int_{0}^{t}\left\vert \eta _{s}\right\vert
^{2}ds+\int_{0}^{t}\int_{E}\left\vert \sigma _{s}\left( e\right) \right\vert
^{2}\lambda \left( de\right) ds \\
&&+\sum_{0\leq t\leq T}\left( \Delta \alpha _{s}\right) ^{2},
\end{eqnarray*}%
$\left( ii\right) $%
\begin{equation*}
\mathbb{E}\left\vert \alpha _{t}\right\vert ^{2}+\mathbb{E}%
\int_{t}^{T}\left\vert \eta _{s}\right\vert ^{2}ds+\mathbb{E}%
\int_{t}^{T}\int_{E}\left\vert \sigma _{s}\left( e\right) \right\vert
^{2}\lambda \left( de\right) ds\leq \mathbb{E}\left\vert \alpha _{T}\right\vert ^{2}+2%
\mathbb{E}\int_{t}^{T}\langle \alpha _{s},\beta _{s}\rangle ds+2\mathbb{E}%
\int_{t}^{T}\langle \alpha _{s},dK_{s}\rangle +\mathbb{E}\int_{t}^{T}\left%
\vert \gamma _{s}\right\vert ^{2}ds.
\end{equation*}%
%
%
%
%
%
%
%
%
%
%
%
%
%
%
%
\end{lemma}

\section{Comparison theorem.}  Given two parameters $\left( \xi ^{1},f^{1},g,T\right) $
and $\left( \xi
^{2},f^{2},g,T\right) ,$ we considere the reflected BDSDEPs, $i=1,2$%
\begin{equation}
\begin{array}{l}
Y_{t}^{i}=\xi
^{i}+\int_{t}^{T}f^{i}(s,Y_{s}^{i},Z_{s}^{i},U_{s}^{i})ds+%
\int_{t}^{T}g(s,Y_{s}^{i},Z_{s}^{i},U_{s}^{i})d\overleftarrow{B}_{s} \\
\\
+\int_{t}^{T}dK_{s}^{i}-\int_{t}^{T}Z_{s}^{i}dW_{s}-\int_{t}^{T}%
\int_{E}U_{s}^{i}\left( e\right) \tilde{\mu}\left( ds,de\right) ,\ 0\leq
t\leq T.%
\end{array}%
\end{equation}

\begin{theorem}
\label{th} Assume that the reflected BDSDEP associated with dates $(\xi
^{1},f^{1},g,T),$ $\left( resp\text{ }\left( \xi ^{^{2}},f^{^{2}},g,T\right)
\right) $ has a solution $(Y_{t}^{1},Z_{t}^{1},K_{t}^{1},U_{t}^{1})_{t\in
\lbrack 0,T)]},$ $($ $resp$ $\left(
Y_{t}^{2},Z_{t}^{2},K_{t}^{2},U_{t}^{2}\right) _{t\in \lbrack 0,T]}).$ Each
one satisfying the assumption \textbf{(H)}, assume moreover that:%
\begin{equation*}
\left\{
\begin{array}{l}
\xi ^{1}\leq \xi ^{2}, \\
\forall t\leq T,\text{ }S_{t}^{1}\leq S_{t}^{2}, \\
f^{1}\left( t,Y_{t},Z_{t},U_{t}\right) \leq f^{^{2}}\left(
t,Y_{t},Z_{t},U_{t}\right) .%
\end{array}%
\right.
\end{equation*}%
Then we have $\mathbb{P}-a.s.,$%
\begin{equation*}
Y_{t}^{1}\leq Y_{t}^{^{2}}.
\end{equation*}%
%
%
%
%
%
%
%
%
%
%
%
%
%
%
%
\end{theorem}

\noindent \textbf{Proof: }Let us show that $\left(
Y_{t}^{1}-Y_{t}^{^{2}}\right) ^{+}=0,$ using the equation (3.1), we get%
\begin{eqnarray*}
\bar{Y}_{t} &=&Y_{t}^{1}-Y_{t}^{2} \\
&=&\bar{\xi}+\int_{t}^{T}\left(
f^{1}(s,Y_{s}^{1},Z_{s}^{1},U_{s}^{1})-f^{2}(s,Y_{s}^{2},Z_{s}^{2},U_{s}^{2})\right) ds+\int_{t}^{T}\left( g(s,Y_{s}^{1},Z_{s}^{1},U_{s}^{1})-g(s,Y_{s}^{2},Z_{s}^{2},U_{s}^{2})\right) d%
\overleftarrow{B}_{s} \\
&&+\int_{t}^{T}\left( dK_{s}^{1}-dK_{s}^{2}\right) -\int_{t}^{T}\bar{Z}%
_{s}dW_{s}-\int_{t}^{T}\int_{E}\bar{U}_{s}\left( e\right) \lambda
\left( ds,de\right) ,
\end{eqnarray*}
where $\bar{\xi}=\xi ^{1}- \xi ^{2}, \ \bar{Z}=Z ^{1}- Z ^{2}$ and
$\bar{U}=U ^{1}- U ^{2}.$\\ Since $\int_{t}^{T}\left(
\bar{Y}_{s}\right) ^{+}\left(
g(s,Y_{s}^{1},Z_{s}^{1},U_{s}^{1})-g(s,Y_{s}^{2},Z_{s}^{2},U_{s}^{2})\right)
d\overleftarrow{B}_{s}$ and $\int_{t}^{T}\left( \bar{Y}_{s}\right) ^{+}\bar{Z%
}_{s}dW_{s}$ are a uniformly integrable martingale then taking expectation,
we get by applying Lemma 2.1%
\begin{eqnarray*}
&&\mathbb{E}\left\vert \left( \bar{Y}_{t}\right) ^{+}\right\vert ^{2}+%
\mathbb{E}\int_{t}^{T}1_{\left\{ \bar{Y}_{s}>0\right\} }\left\vert
\left\vert \bar{Z}_{s}\right\vert \right\vert ^{2}ds+\mathbb{E}%
\int_{t}^{T}\int_{E}1_{\left\{ \bar{Y}_{s}>0\right\} }\left\vert \bar{U}%
_{s}\left( e\right) \right\vert ^{2}\lambda \left( de\right) ds \\
&\leq &\mathbb{E}\left\vert \left( \bar{\xi}\right) ^{+}\right\vert ^{2}+2%
\mathbb{E}\int_{t}^{T}\left( \bar{Y}_{s}\right) ^{+}\left(
f^{1}(s,Y_{s}^{1},Z_{s}^{1},U_{s}^{1})-f^{2}(s,Y_{s}^{2},Z_{s}^{2},U_{s}^{2})\right) ds
\\
&&+2\mathbb{E}\int_{t}^{T}\left( \bar{Y}_{s}\right) ^{+}\left(
dK_{s}^{1}-dK_{s}^{2}\right) +\mathbb{E}\int_{t}^{T}1_{_{\left\{ \bar{Y}%
_{s}>0\right\} }}\left\vert \left\vert
g(s,Y_{s}^{1},Z_{s}^{1},U_{s}^{1})-g(s,Y_{s}^{2},Z_{s}^{2},U_{s}^{2})\right%
\vert \right\vert ^{2}ds.
\end{eqnarray*}
Since%
\begin{equation*}
\left\{
\begin{array}{l}
\left( \xi ^{1}-\xi ^{2}\right) ^{+}=0, \\
\\
\int_{t}^{T}\left( \bar{Y}_{s}\right) ^{+}\left(
dK_{s}^{1}-dK_{s}^{2}\right) =-\int_{t}^{T}\left( Y_{s}^{1}-Y_{s}^{2}\right)
^{+}dK_{s}^{2}\leq 0,%
\end{array}%
\right.
\end{equation*}%
we get%
\begin{eqnarray*}
&&\mathbb{E}\left\{ \left\vert \left( \bar{Y}_{t}\right) ^{+}\right\vert
^{2}+\int_{t}^{T}1_{\left\{ \bar{Y}_{s}>0\right\} }\left\vert \left\vert
\bar{Z}_{s}\right\vert \right\vert ^{2}ds+\int_{t}^{T}\int_{E}1_{\left\{
\bar{Y}_{s}>0\right\} }\left\vert \bar{U}_{s}\left( e\right) \right\vert ^{2}%
\tilde{\mu}\left( de\right) ds\right\}  \\
&\leq &2\mathbb{E}\int_{t}^{T}\left( \bar{Y}_{s}\right) ^{+}\left(
f^{1}(s,Y_{s}^{1},Z_{s}^{1},U_{s}^{1})-f^{2}(s,Y_{s}^{2},Z_{s}^{2},U_{s}^{2})\right) ds
\\
&&+\mathbb{E}\int_{t}^{T}1_{_{\left\{ \bar{Y}_{s}>0\right\} }}\left\vert
\left\vert
g(s,Y_{s}^{1},Z_{s}^{1},U_{s}^{1})-g(s,Y_{s}^{2},Z_{s}^{2},U_{s}^{2})\right%
\vert \right\vert ^{2}ds,
\end{eqnarray*}%
we obtain, by hypothesis $\left( H.2\right) ,$ and Young's inequality the
following inequality%
\begin{eqnarray*}
&&2\mathbb{E}\int_{t}^{T}\left( \bar{Y}_{s}\right) ^{+}\left(
f^{1}(s,Y_{s}^{1},Z_{s}^{1},U_{s}^{1})-f^{2}(s,Y_{s}^{2},Z_{s}^{2},U_{s}^{2})\right) ds
\\
&\leq &(2C+2C^2\epsilon) \mathbb{E}\int_{t}^{T}\left\vert \bar{Y}_{s}^{+}\right\vert
^{2}ds+\epsilon ^{-1}\mathbb{E}\int_{t}^{T}\left( \left\vert \bar{Z}_{s}\right\vert
^{2}+\int_{E}\left\vert \bar{U}_{s}\right\vert ^{2}\lambda \left( de\right)
\right) ds,
\end{eqnarray*}%
also we applying the assumption $\left( H.2\right) $ for $g,$ we get%
\begin{equation*}
\left\vert \left\vert
g(s,Y_{s}^{1},Z_{s}^{1},U_{s}^{1})-g(s,Y_{s}^{2},Z_{s}^{2},U_{s}^{2})\right%
\vert \right\vert ^{2}\leq C\left\vert \bar{Y}_{s}\right\vert
^{2}ds+\alpha
\left\{ \left\vert \bar{Z}_{s}\right\vert ^{2}+\left\vert \bar{U}%
_{s}\right\vert ^{2}\right\} .
\end{equation*}%
Then, we have the following inequality%
\begin{eqnarray*}
&&\mathbb{E}\left\{ \left\vert \left( \bar{Y}_{t}\right) ^{+}\right\vert
^{2}+\int_{t}^{T}1_{\left\{ \bar{Y}_{s}>0\right\} }\left\vert \left\vert
\bar{Z}_{t}\right\vert \right\vert ^{2}ds+\int_{t}^{T}\int_{E}1_{\left\{
\bar{Y}_{s}>0\right\} }\left\vert \bar{U}_{s}\left( e\right) \right\vert ^{2}%
\tilde{\mu}\left( de\right) ds\right\}  \\
&\leq &(2C+2C^2\epsilon) \mathbb{E}\int_{t}^{T}\left\vert \bar{Y}_{s}^{+}\right\vert
^{2}ds+\epsilon ^{-1}\mathbb{E}\int_{t}^{T}\left( \left\vert \bar{Z}_{s}\right\vert
^{2}+\int_{E}\left\vert \bar{U}_{s}\right\vert ^{2}\lambda \left( de\right)
\right) ds \\
&&+C\mathbb{E}\int_{t}^{T}1_{_{\left\{ \bar{Y}_{s}>0\right\}
}}\left\vert \bar{Y}_{s}\right\vert ^{2}ds+\alpha \mathbb{E}
\int_{t}^{T}\left\{ 1_{\left\{ \bar{Y}_{s}>0\right\} }\left\vert \bar{Z}%
_{s}\right\vert ^{2}+\int_{E}1_{_{\left\{ \bar{Y}_{s}>0\right\} }}\left\vert
\bar{U}_{s}\right\vert ^{2}\lambda \left( de\right) \right\} ds, \\
&\leq&\left(2C^2 \epsilon +3C\right) \mathbb{E}\int_{t}^{T}\left%
\vert \bar{Y}_{s}^{+}\right\vert ^{2}ds+\left( \epsilon
^{-1}+\alpha \right) \mathbb{E}\left\{ \int_{t}^{T}1_{_{\left\{
\bar{Y}_{s}>0\right\} }}\left\vert \bar{Z}_{s}\right\vert
^{2}ds+\int_{t}^{T}\int_{E}1_{_{\left\{ \bar{Y}_{s}>0\right\}
}}\left\vert \bar{U}_{s}\right\vert ^{2}\lambda \left( de\right)
ds\right\} ,
\end{eqnarray*}%
choosing $\epsilon $ such that
$0<\epsilon^{-1}+\alpha\leq1$, we have
\begin{equation*}
\mathbb{E}\left\vert \left( \bar{Y}_{t}\right) ^{+}\right\vert
^{2}\leq \left( 3C+2C^2\epsilon\right)
\mathbb{E}\int_{t}^{T}\left\vert \bar{Y}_{s}^{+}\right\vert ^{2}ds,
\end{equation*}%
using Gronwall's lemma implies that%
\begin{equation*}
\mathbb{E}\left[ \left\vert \left( \bar{Y}_{t}\right) ^{+}\right\vert ^{2}%
\right] =0,
\end{equation*}%
finally, we have%
\begin{equation*}
Y_{t}^{1}\leq Y_{t}^{^{2}}.
\end{equation*}%
\eop

\section{Reflected BDSDEPs with continuous coefficient.}

\qquad In this section we are interested in weakening the conditions on $f$.
We assume that $f$ and $g$ satisfy the following assumptions:

\noindent \vskip0.1cm\noindent \textbf{(H.6) }There exists $0<\alpha<1$ and $C>0$ $s.t$. for
all $\left( t,\omega ,y,z,u\right) $ $\in \left[ 0,T\right] \times \Omega
\times
\mathbb{R}
\times
\mathbb{R}
^{d}\times L^{2}\left( E,\mathcal{E},\lambda ,%
\mathbb{R}
\right) ,$%
\begin{equation*}
\left\{
\begin{array}{l}
\left\vert f(t,\omega ,y,z,u)\right\vert \leq C\left( 1+\left\vert
y\right\vert +\left\vert z\right\vert +\left\vert u\right\vert \right) , \\
\left\vert g(t,\omega ,y,z,u)-g(t,\omega ,y^{^{\prime }},z^{^{\prime
}},u^{^{\prime }})\right\vert ^{2}\leq C\left\vert y-y^{^{\prime
}}\right\vert ^{2}+\alpha \left\{ \left\vert z-z^{^{\prime }}\right\vert
^{2}+\left\vert u-u^{^{\prime }}\right\vert ^{2}\right\} .%
\end{array}%
\right.
\end{equation*}%
\noindent \vskip0.1cm\noindent \textbf{(H.7) }For fixed $\omega $ and $t$, $%
f(t,\omega ,\cdot ,\cdot ,\cdot )$ is continuous.

The theree next Lemmas will be useful in the sequel.

we recall the following classical lemma. It can be
proved by adapting the proof given in J. J. Alibert and K. Bahlali [1].
\begin{lemma}
\label{Lemma 2} Let $f:\Omega \times \lbrack 0,T]\times
\mathbb{R}
\times
\mathbb{R}
^{d}\times \mathcal{L}^{2}\left( E,\mathcal{E},\lambda ,%
\mathbb{R}
\right) \rightarrow
\mathbb{R}
$ be a mesurable function such that:

\begin{enumerate}
\item For a.s. every $\left( t,\omega \right) \in \lbrack 0,T]\times \Omega
, $ $f\left( t,\omega ,y,z,u\right) $ is a continuous.

\item There exists a constant $C> 0$ such that for every $\left( t,\omega
,y,z,u\right) \in \lbrack 0,T]\times \Omega \times
\mathbb{R}
\times
\mathbb{R}
^{d}\times L^{2}\left( E,\mathcal{E},\lambda ,%
\mathbb{R}
\right) ,$ $\left\vert f\left( t,\omega ,y,z,u\right) \right\vert \leq
C\left( 1+\left\vert y\right\vert +\left\vert z\right\vert +\left\vert
u\right\vert \right) .$
\end{enumerate}

Then exists the sequence of fonction $f_{n}$%
\begin{equation*}
f_{n}\left( t,\omega ,y,z,u\right) =\inf_{\left( y^{^{\prime }},z^{\prime
},u^{^{\prime }}\right) \in \mathcal{B}^{2}\left(
\mathbb{R}
\right) }\left[ f\left( t,\omega ,y^{^{\prime }},z^{\prime },u^{^{\prime
}}\right) +n\left( \left\vert y-y^{^{\prime }}\right\vert +\left\vert
z-z^{\prime }\right\vert +\left\vert u-u^{^{\prime }}\right\vert \right) %
\right] ,
\end{equation*}%
is well defined for each $n\geq C$, and it satisfies, $d\mathbb{P}\times
dt-a.s.$

$\left( i\right) $ Linear growth: $\forall n\geq 1,$ $\left( y,z,u\right)
\in \mathbb{R}\times\mathbb{R}^d\times\mathcal{L}^2 $, $\left\vert f_{n}\left( t,\omega ,y,z,u\right) \right\vert \leq
C\left( 1+\left\vert y\right\vert +\left\vert z\right\vert +\left\vert
u\right\vert \right) .$

$\left( ii\right) $ Monotonicity in $n:\forall $ $y,z,u,$  $f_{n}\left(
t,\omega ,y,z,u\right) $ is increases in $n.$

$\left( iii\right) $ Convergence: $\forall \left( t,\omega ,y,z,u\right) \in %
\left[ 0,T\right] \times \Omega \times \mathcal{B}^{2}\left(
\mathbb{R}
\right) $, if $\left( t,\omega ,y_{n},z_{n},u_{n}\right) \rightarrow \left(
t,\omega ,y,z,u\right) ,$ then $f_{n}\left( t,\omega
,y_{n},z_{n},u_{n}\right) \rightarrow f\left( t,\omega ,y,z,u\right) .$

$\left( iv\right) $ Lipschitz condition: $\forall n\geq 1,$ $\left( t,\omega
\right) \in \lbrack 0,T]\times \Omega ,$ $\forall \left( y,z,u\right) \in
\mathcal{B}^{2}\left(
\mathbb{R}
\right) $ and $\left( y^{^{\prime }},z^{\prime },u^{^{\prime }}\right) \in
\mathcal{B}^{2}\left(
\mathbb{R}
\right) $, we have%
\begin{equation*}
\left\vert f_{n}\left( t,\omega ,y,z,u\right) -f_{n}\left( t,\omega
,y^{^{\prime }},z^{\prime },u^{^{\prime }}\right) \right\vert \leq
n\left\vert y-y^{^{\prime }}\right\vert +\left\vert z-z^{\prime }\right\vert
+\left\vert u-u^{^{\prime }}\right\vert .
\end{equation*}%
%
%
%
%
%
%
%
%
%
%
%
%
%
%
%
\end{lemma}

Now given $\xi \in $ $\mathbb{L}^{2}$, $n\in N$, we consider $\left(
Y^{n},Z^{n},K^{n},U^{n}\right) $ and $\left( \text{resp }\left(
V,N,K,M\right) \right) $ be solutions of the following reflected BDSDEPs:%
\begin{equation}
\left\{
\begin{array}{l}
Y_{t}^{n}=\xi
+\int_{t}^{T}f_{n}(s,Y_{s}^{n},Z_{s}^{n},U_{s}^{n})ds+%
\int_{t}^{T}g(s,Y_{s}^{n},Z_{s}^{n},U_{s}^{n})d\overleftarrow{B}_{s} \\
\\
+\int_{t}^{T}dK_{s}^{n}-\int_{t}^{T}Z_{s}^{n}dW_{s}-\int_{t}^{T}%
\int_{E}U_{s}^{n}\left( e\right) \tilde{\mu}\left( ds,de\right) ,\ 0\leq
t\leq T, \\
\\
S_{t}\leq Y_{t}^{n},\text{ }0\leq t\leq T,\text{\quad and\quad }%
\int_{0}^{T}\left( Y_{t}^{n}-S_{t}\right) dK_{t}^{n}=0.%
\end{array}%
\right.
\end{equation}
respectively
\begin{equation*}
\left\{
\begin{array}{l}
V_{t}=\xi
+\int_{t}^{T}H(s,V_{s},N_{s},M_{s})ds+\int_{t}^{T}g(s,V_{s},N_{s},M_{s})d%
\overleftarrow{B}_{s} \\
\\
+\int_{t}^{T}dK_{s}-\int_{t}^{T}N_{s}dW_{s}-\int_{t}^{T}\int_{E}M_{s}\left(
e\right) \tilde{\mu}\left( ds,de\right) ,\ 0\leq t\leq T, \\
\\
S_{t}\leq V_{t},\text{ }0\leq t\leq T,\text{\quad and\quad }%
\int_{0}^{T}\left( V_{t}-S_{t}\right) dK_{t}=0,%
\end{array}%
\right.
\end{equation*}%
where $H(s,\omega ,V,N,M)=C\left( 1+\left\vert V\right\vert
+\left\vert N\right\vert +\left\vert M\right\vert \right) .$

\begin{lemma}
\label{Lemma 3} $\left( i\right) $ $a.s.$ for all, $t$ and $\forall
n\leq m,$ $Y_{t}^{n}\leq Y_{t}^{m}\leq V_{t}.$

$\left( ii\right) $ Assume that $\left( H.1\right) ,$ $\left( H.3\right)
-\left( H.7\right) $ is in force. Then there exists a constant $A>0$
depending only on $C,$ $\alpha ,$ $\xi $ and $T$ such that:%
\begin{equation*}
\left\vert \left\vert U^{n}\right\vert \right\vert _{\mathcal{L}^{2}\left(
0,T,\tilde{\mu},%
\mathbb{R}
\right) }\leq A,\text{\qquad \qquad \qquad }\left\vert \left\vert
Z^{n}\right\vert \right\vert _{\mathcal{M}^{2}\left( 0,T,%
\mathbb{R}
^{d}\right) }\leq A.
\end{equation*}%
%
%
%
%
%
%
%
%
\end{lemma}

\noindent \textbf{Proof: }The prove of the\textbf{\ }$\left( i\right) $
follow from comparison theorem. It remains to prove $\left( ii\right) $, by
lemma 2.1, we have%
\begin{eqnarray*}
&&\mathbb{E}\left\vert Y_{t}^{n}\right\vert ^{2}+\mathbb{E}%
\int_{t}^{T}\left\vert Z_{s}^{n}\right\vert ^{2}ds+\mathbb{E}%
\int_{t}^{T}\int_{E}\left\vert U_{s}^{n}\left( e\right) \right\vert
^{2}\lambda \left( de\right) ds \\
&\leq &\mathbb{E}\left\vert \xi \right\vert ^{2}+2\mathbb{E}%
\int_{t}^{T}Y_{s}^{n}f_{n}(s,Y_{s}^{n},Z_{s}^{n},U_{s}^{n})ds+2\mathbb{E}%
\int_{t}^{T}Y_{s}^{n}dK_{s}^{n}+\mathbb{E}\int_{t}^{T}\left\vert \left\vert
g(s,Y_{s}^{n},Z_{s}^{n},U_{s}^{n})\right\vert \right\vert ^{2}ds.  \notag
\end{eqnarray*}%
By $\left( i\right) $ in lemma 4.1, we have%
\begin{eqnarray*}
2\mathbb{E}\int_{t}^{T}Y_{s}^{n}f_{n}(s,Y_{s}^{n},Z_{s}^{n},U_{s}^{n})ds
&\leq &2C\mathbb{E}\int_{t}^{T}Y_{s}^{n}\left( 1+\left\vert
Y_{s}^{n}\right\vert +\left\vert Z_{s}^{n}\right\vert +\left\vert
U_{s}^{n}\right\vert \right) ds \\
&\leq &\mathbb{E}\int_{t}^{T}\left\vert Y_{s}^{n}\right\vert ^{2}ds+TC^{2}+2C%
\mathbb{E}\int_{t}^{T}\left\vert Y_{s}^{n}\right\vert ^{2}ds+\frac{C^{2}}{%
\gamma _{1}}\mathbb{E}\int_{t}^{T}\left\vert Y_{s}^{n}\right\vert ^{2}ds \\
&&+\gamma _{1}\mathbb{E}\int_{t}^{T}\left\vert Z_{s}^{n}\right\vert ^{2}ds+%
\frac{C^{2}}{\gamma _{2}}\mathbb{E}\int_{t}^{T}\left\vert
Y_{s}^{n}\right\vert ^{2}ds+\gamma _{2}\mathbb{E}\int_{t}^{T}\int_{E}\left%
\vert U_{s}^{n}\left( e\right) \right\vert ^{2}\lambda \left( de\right) ds,
\\
&\leq &\left( 1+2C+\frac{C^{2}}{\gamma _{1}}+\frac{C^{2}}{\gamma _{2}}%
\right) \mathbb{E}\int_{t}^{T}\left\vert Y_{s}^{n}\right\vert ^{2}ds+TC^{2}
\\
&&+\gamma _{1}\mathbb{E}\int_{t}^{T}\left\vert Z_{s}^{n}\right\vert
^{2}ds+\gamma _{2}\mathbb{E}\int_{t}^{T}\int_{E}\left\vert U_{s}^{n}\left(
e\right) \right\vert ^{2}\lambda \left( de\right) ds,
\end{eqnarray*}%
also by the hypothesis associated with $g$, we get%
\begin{eqnarray*}
\left\vert \left\vert g(s,Y_{s}^{n},Z_{s}^{n},U_{s}^{n})\right\vert
\right\vert ^{2} &\leq &\left( 1+\epsilon \right) \left\vert \left\vert
g(s,Y_{s}^{n},Z_{s}^{n},U_{s}^{n})-g(s,0,0,0)\right\vert \right\vert ^{2}+%
\frac{1+\epsilon }{\epsilon }\left\vert \left\vert g(s,0,0,0)\right\vert
\right\vert ^{2}, \\
&\leq &\left( 1+\epsilon \right) C\left\vert Y_{s}^{n}\right\vert
^{2}+\left( 1+\epsilon \right) \alpha \left\{ \left\vert
Z_{s}^{n}\right\vert ^{2}+\left\vert U_{s}^{n}\right\vert ^{2}\right\} +%
\frac{1+\epsilon }{\epsilon }\left\vert \left\vert g(s,0,0,0)\right\vert
\right\vert ^{2}.
\end{eqnarray*}%
Chossing $\gamma _{1}=\gamma _{2}=\frac{\epsilon ^{2}}{2}$. Then, we obtain
the following inequality%
\begin{equation*}
\begin{array}{l}
\mathbb{E}\left( \left\vert Y_{t}^{n}\right\vert ^{2}+\int_{t}^{T}\left\vert
Z_{s}^{n}\right\vert ^{2}ds+\int_{t}^{T}\int_{E}\left\vert U_{s}^{n}\left(
e\right) \right\vert ^{2}\lambda \left( de\right) ds\right)  \\
\leq \mathbb{E}\left\vert \xi \right\vert ^{2}+TC^{2}+\left( 1+2C+\frac{%
4C^{2}}{\epsilon ^{2}}+\left( 1+\epsilon \right) C\right) \mathbb{E}%
\int_{0}^{T}\left\vert Y_{s}^{n}\right\vert
^{2}ds+2\int_{0}^{T}Y_{s}^{n}dK_{s}^{n} \\
+\left( \frac{\epsilon ^{2}}{2}+\left( 1+\epsilon \right) \alpha \right)
\left\{ \mathbb{E}\int_{t}^{T}\left\vert Z_{s}^{n}\right\vert ^{2}ds+\mathbb{%
E}\int_{0}^{T}\int_{E}\left\vert U_{s}^{n}\left( e\right) \right\vert
^{2}\lambda \left( de\right) ds\right\} +\frac{1+\epsilon }{\epsilon }%
\mathbb{E}\int_{0}^{T}\left\vert \left\vert g(s,0,0,0)\right\vert
\right\vert ^{2}ds.%
\end{array}%
\end{equation*}%
Consequently, we have%
\begin{eqnarray*}
&&\mathbb{E}\int_{t}^{T}\left( \left\vert Z_{s}^{n}\right\vert
^{2}+\int_{E}\left\vert U_{s}^{n}\left( e\right) \right\vert ^{2}\lambda
\left( de\right) \right) ds  \notag \\
&\leq &\left( \frac{\epsilon ^{2}}{2}+\left( 1+\epsilon \right) \alpha
\right) \mathbb{E}\left\{ \int_{t}^{T}\left\vert Z_{s}^{n}\right\vert
^{2}ds+\int_{t}^{T}\int_{E}\left\vert U_{s}^{n}\left( e\right) \right\vert
^{2}\lambda \left( de\right) ds\right\} +\Lambda +\theta \mathbb{E}%
\left\vert K_{T}^{n}-K_{t}^{n}\right\vert ^{2},
\end{eqnarray*}
where

$\Lambda =\mathbb{E}\left\vert \xi \right\vert ^{2}+TC^{2}+\frac{1+\epsilon
}{\epsilon }\mathbb{E}\int_{t}^{T}\left\vert \left\vert
g(s,0,0,0)\right\vert \right\vert ^{2}ds+\frac{1}{\theta }\mathbb{E}\left(
\sup_{0\leq s\leq T}\left( S_{s}\right) ^{2}\right) +T\left( 1+2C+\frac{%
4C^{2}}{\epsilon ^{2}}+\left( 1+\epsilon \right) C\right) \mathbb{E}\left(
\sup_{t}\left\vert Y_{t}^{n}\right\vert ^{2}\right).$ Now chossing $%
\epsilon $ and $\alpha $ such that $0\leq \frac{\epsilon ^{2}}{2}+\left(
1+\epsilon \right) \alpha <1$, we obtain%
\begin{equation}
\mathbb{E}\int_{t}^{T}\left\vert Z_{s}^{n}\right\vert ^{2}ds+\mathbb{E}%
\int_{t}^{T}\int_{E}\left\vert U_{s}^{n}\left( e\right) \right\vert
^{2}\lambda \left( de\right) ds\leq \Lambda +\theta \mathbb{E}\left\vert
K_{T}^{n}-K_{t}^{n}\right\vert ^{2}.
\end{equation}
On the other hand, we have from Eq.$\left( 4.1\right) $%
\begin{eqnarray*}
K_{T}^{n}-K_{t}^{n} &=&Y_{t}^{n}-\xi
-\int_{t}^{T}f_{n}(s,Y_{s}^{n},Z_{s}^{n},U_{s}^{n})ds-%
\int_{t}^{T}g(s,Y_{s}^{n},Z_{s}^{n},U_{s}^{n})d\overleftarrow{B}_{s} \\
&&+\int_{t}^{T}Z_{s}^{n}dW_{s}+\int_{t}^{T}\int_{E}U_{s}^{n}\left( e\right)
\tilde{\mu}\left( ds,de\right) .
\end{eqnarray*}%
Using the Hölder's inequality and assupmtion $\left( H.6\right) $, we have%
\begin{equation*}
\mathbb{E}\left\vert K_{T}^{n}-K_{t}^{n}\right\vert ^{2}\leq
C_{1}+C_{2}\left( \mathbb{E}\int_{t}^{T}\left\vert Z_{s}^{n}\right\vert
^{2}ds+\mathbb{E}\int_{t}^{T}\int_{E}\left\vert U_{t}^{n}\right\vert
^{2}\lambda \left( de\right) ds\right) ,
\end{equation*}%
From inequality $\left( 4.2\right) $, we get%
\begin{equation*}
\mathbb{E}\int_{0}^{T}\left( \left\vert Z_{s}^{n}\right\vert
^{2}+\int_{E}\left\vert U_{s}^{n}\left( e\right) \right\vert ^{2}\lambda
\left( de\right) \right) ds\leq \Lambda +\theta C_{1}+\theta C_{2}\mathbb{E}%
\int_{t}^{T}\left( \left\vert Z_{s}^{n}\right\vert ^{2}+\int_{E}\left\vert
U_{t}^{n}\right\vert ^{2}\lambda \left( de\right) \right) ds,
\end{equation*}%
Finally chossing $\theta $ such that $0\leq \theta C_{2}\leq 1$, we obtain%
\begin{equation*}
\mathbb{E}\int_{t}^{T}\left\vert Z_{s}^{n}\right\vert ^{2}ds+\mathbb{E}%
\int_{t}^{T}\int_{E}\left\vert U_{s}^{n}\left( e\right) \right\vert
^{2}\lambda \left( de\right) ds\leq \Lambda +\theta C_{1}<\infty .
\end{equation*}%
The prove of lemma 4.2 is complet.\eop

\begin{lemma}
\label{Lemma 4} Assume that $\left( H.1\right) ,$ $\left( H.3\right) -\left(
H.7\right) $ is in force. Then the sequence $\left( Z^{n},U^{n}\right) $
converges a.s. in $\mathcal{M}^{2}\left( 0,T,%
\mathbb{R}
^{d}\right) \times \mathcal{L}^{2}\left( 0,T,\tilde{\mu},%
\mathbb{R}
\right) .$

\end{lemma}

\noindent \textbf{Proof: }Let $n_{0}\geq C$. From Eq.$\left(
4.1\right) $,
we deduce that there exists a process $Y\in \mathcal{S}^{2}\left( 0,T,%
\mathbb{R}
\right) $ such that $Y^{n}\rightarrow Y$ a.s., as $n\rightarrow \infty $.
Applying Lemma 2.1 to $\left\vert Y_{t}^{n}-Y_{t}^{m}\right\vert ^{2}$, for $%
n,m\geq n_{0}$%
\begin{equation*}
\begin{array}{l}
\mathbb{E}\left( \left\vert Y_{t}^{n}-Y_{t}^{m}\right\vert
^{2}+\int_{t}^{T}\left\vert Z_{s}^{n}-Z_{s}^{m}\right\vert
^{2}ds+\int_{t}^{T}\int_{E}\left\vert U_{s}^{n}\left( e\right)
-U_{s}^{m}\left( e\right) \right\vert ^{2}\lambda \left( de\right) ds\right)
\\
\leq 2\mathbb{E}\int_{t}^{T}\left( Y_{s}^{n}-Y_{s}^{m}\right) \left(
f_{n}(s,Y_{s}^{n},Z_{s}^{n},U_{s}^{n})-f_{m}(s,Y_{s}^{m},Z_{s}^{m},U_{s}^{m})\right) ds
\\
+2\mathbb{E}\int_{t}^{T}\left( Y_{s}^{n}-Y_{s}^{m}\right) \left(
dK_{s}^{n}-dK_{s}^{m}\right) +\mathbb{E}\int_{t}^{T}\left\vert \left\vert
g(s,Y_{s}^{n},Z_{s}^{n},U_{s}^{n})-g(s,Y_{s}^{m},Z_{s}^{m},U_{s}^{m})\right%
\vert \right\vert ^{2}ds.%
\end{array}%
\end{equation*}%
Since $\int_{t}^{T}\left( Y_{s}^{n}-Y_{s}^{m}\right) \left(
dK_{s}^{n}-dK_{s}^{m}\right) \leq 0$, we deduce that%
\begin{eqnarray*}
&&\mathbb{E}\int_{t}^{T}\left\vert Z_{t}^{n}-Z_{t}^{m}\right\vert ^{2}ds+%
\mathbb{E}\int_{t}^{T}\int_{E}\left\vert U_{s}^{n}\left( e\right)
-U_{s}^{m}\left( e\right) \right\vert ^{2}\lambda \left( de\right) ds \\
&\leq &2\mathbb{E}\int_{t}^{T}\left( Y_{s}^{n}-Y_{s}^{m}\right) \left(
f_{n}(s,Y_{s}^{n},Z_{s}^{n},U_{s}^{n})-f_{m}(s,Y_{s}^{m},Z_{s}^{m},U_{s}^{m})\right) ds
\\
&&+\mathbb{E}\int_{t}^{T}\left\vert \left\vert
g(s,Y_{s}^{n},Z_{s}^{n},U_{s}^{n})-g(s,Y_{s}^{m},Z_{s}^{m},U_{s}^{m})\right%
\vert \right\vert ^{2}ds.
\end{eqnarray*}%
Using Hölder's inequality and assumption $\left( H.6\right) $ for $g$, we
deduce that%
\begin{eqnarray*}
&&\left( 1-\alpha \right) \mathbb{E}\left\{ \int_{t}^{T}\left\vert
Z_{t}^{n}-Z_{t}^{m}\right\vert ^{2}ds+\int_{t}^{T}\int_{E}\left\vert
U_{s}^{n}\left( e\right) -U_{s}^{m}\left( e\right) \right\vert ^{2}\lambda
\left( de\right) ds\right\}  \\
&\leq &2\mathbb{E}\left( \int_{t}^{T}\left\vert
f_{n}(s,Y_{s}^{n},Z_{s}^{n},U_{s}^{n})-f_{m}(s,Y_{s}^{m},Z_{s}^{m},U_{s}^{m})\right\vert ^{2}ds\right) ^{%
\frac{1}{2}}\mathbb{E}\left( \int_{t}^{T}\left\vert
Y_{s}^{n}-Y_{s}^{m}\right\vert ^{2}ds\right) ^{\frac{1}{2}}+C\mathbb{E}%
\int_{t}^{T}\left\vert Y_{s}^{n}-Y_{s}^{m}\right\vert ^{2}ds.
\end{eqnarray*}%
Applying assumption $\left( H.6\right) $ for $f$ and the boundedness of the
sequence $\left( Y^{n},Z^{n},U^{n}\right) $, we deduce that%
\begin{equation*}
\left( 1-\alpha \right) \left\{ \mathbb{E}\int_{t}^{T}\left\vert
Z_{t}^{n}-Z_{t}^{m}\right\vert ^{2}ds+\mathbb{E}\int_{t}^{T}\int_{E}\left%
\vert U_{s}^{n}\left( e\right) -U_{s}^{m}\left( e\right) \right\vert
^{2}\lambda \left( de\right) ds\right\} \leq C^{te}\mathbb{E}%
\int_{t}^{T}\left\vert Y_{s}^{n}-Y_{s}^{m}\right\vert ^{2}ds,
\end{equation*}%
where the constant $C^{te}>0$ depend only $\xi,\ C,$ $\alpha $ and $T.$

Which yields that $\left( Z^{n}\right) _{n\geq 0}$ respectively $\left(
U^{n}\right) _{n\geq 0}$ is a Cauchy sequence in $\mathcal{M}^{2}\left( 0,T,%
\mathbb{R}
^{d}\right) ,$ respectively in $\mathcal{L}^{2}\left( 0,T,\tilde{\mu},%
\mathbb{R}
\right) .$ Then there exists $\left( Z,U\right) \in \mathcal{M}^{2}\left(
0,T,%
\mathbb{R}
^{d}\right) \times \mathcal{L}^{2}\left( 0,T,\tilde{\mu},%
\mathbb{R}
\right) $ such that%
\begin{equation*}
\mathbb{E}\int_{0}^{T}\left\vert Z_{s}^{n}-Z_{s}\right\vert ^{2}ds+\mathbb{E}%
\int_{0}^{T}\int_{E}\left\vert U_{s}^{n}\left( e\right) -U_{s}\left(
e\right) \right\vert ^{2}\lambda \left( de\right) ds\rightarrow 0,\text{%
\quad }as\text{ }n\rightarrow \infty .
\end{equation*}%
\eop

\begin{theorem}
\label{th1}Assume that $\left( H.1\right) ,$ $\left( H.3\right) -\left(
H.7\right) $ holds. Then Eq $\left( 1.1\right) $ admits a solution $\left(
Y,Z,U,K\right) \in \mathcal{D}^{2}\left(
\mathbb{R}
\right) $. Moreover there is a minimal solution $\left( Y^{\ast },Z^{\ast
},U^{\ast }\right) $ of RBDSDEP $\left( 1.1\right) $ in the sense that for
any other solution $\left( Y,Z,U\right) $ of Eq. $\left( 1.1\right) $, we
have $Y^{\ast }\leq Y$.
\end{theorem}


\noindent \textbf{Proof :}

From Eq.$\left( 4.1\right) $, it's readily seen that $\left(
Y^{n}\right) $
converges in $\mathcal{S}^{2}\left( 0,T,%
\mathbb{R}
\right) $, $dt\otimes d\mathbb{P}-a.s.$ to $Y\in \mathcal{S}^{2}\left( 0,T,%
\mathbb{R}
\right) $. Otherwise thanks to Lemma 4.3 there exists two
subsequences still
noted as the whole sequence $\left( Z^{n}\right) _{n\geq 0}$ respectively $%
\left( U^{n}\right) _{n\geq 0}$ such that%
\begin{equation*}
\begin{tabular}{lll}
$\mathbb{E}\int_{0}^{T}\left\vert Z_{s}^{n}-Z_{s}\right\vert
^{2}ds\rightarrow 0\text{ as }n\rightarrow \infty ,$ & and & $\mathbb{E}%
\int_{0}^{T}\int_{E}\left\vert U_{s}^{n}\left( e\right) -U_{s}\left(
e\right) \right\vert ^{2}\lambda \left( de\right) ds\rightarrow 0,$ $\text{%
as }n\rightarrow \infty .$%
\end{tabular}%
\end{equation*}%
Applying Lemma 4.1, we have $f_{n}\left( t,Y^{n},Z^{n},U^{n}\right)
\rightarrow f\left( t,Y,Z,U\right) $ and the linear growth of $f_{n}$ implies%
\begin{equation*}
\left\vert f_{n}\left( t,Y_{t}^{n},Z_{t}^{n},U_{t}^{n}\right) \right\vert
\leq C\left( 1+\sup_{n}\left( \left\vert Y_{t}^{n}\right\vert
+|Z_{t}^{n}|+|U_{t}^{n}|\right) \right) \in \mathbb{L}^{1}\left( \left[ 0,T%
\right] ;dt\right) .
\end{equation*}%
Thus by Lebesgue's dominated convergence theorem, we deduce that for almost
all $\omega $ and uniformly in $t,$ we have%
\begin{equation*}
\mathbb{E}\int_{t}^{T}f_{n}\left( s,Y_{s}^{n},Z_{s}^{n},U_{s}^{n}\right)
ds\rightarrow \mathbb{E}\int_{t}^{T}f\left( s,Y_{s},Z_{s},U_{s}\right) ds.
\end{equation*}%
We have by $\left( H.6\right) $ the following estimation%
\begin{eqnarray*}
&&\mathbb{E}\int_{t}^{T}\left\vert \left\vert
g(s,Y_{s}^{n},Z_{s}^{n},U_{s}^{n})-g(s,Y_{s},Z_{s},U_{s})\right\vert
\right\vert ^{2}ds \\
&\leq &C\mathbb{E}\int_{t}^{T}\left\vert Y_{s}^{n}-Y_{s}\right\vert
^{2}ds+\alpha \mathbb{E}\int_{t}^{T}\left\vert Z_{s}^{n}-Z_{s}\right\vert
^{2}ds+\alpha \mathbb{E}\int_{t}^{T}\int_{E}\left\vert U_{s}^{n}\left(
e\right) -U_{s}\left( e\right) \right\vert ^{2}\lambda \left( de\right)
ds\rightarrow 0,\text{ as }n\rightarrow \infty ,
\end{eqnarray*}%
using Burkholder-Davis-Gundy inequality, we have%
\begin{equation*}
\left\{
\begin{array}{l}
\mathbb{E}\sup_{0\leq t\leq T}\left\vert
\int_{t}^{T}Z_{s}^{n}dW_{s}-\int_{t}^{T}Z_{s}dW_{s}\right\vert
^{2}\rightarrow 0, \\
\\
\mathbb{E}\sup_{0\leq t\leq T}\left\vert \int_{t}^{T}\int_{E}U_{s}^{n}\left(
e\right) \tilde{\mu}\left( ds,de\right) -\int_{t}^{T}\int_{E}U_{s}\left(
e\right) \tilde{\mu}\left( ds,de\right) \right\vert ^{2}\rightarrow 0, \\
\\
\mathbb{E}\sup_{0\leq t\leq T}\left\vert
\int_{t}^{T}g(s,Y_{s}^{n},Z_{s}^{n},U_{s}^{n})d\overleftarrow{B}%
_{s}-\int_{t}^{T}g(s,Y_{s},Z_{s},U_{s})d\overleftarrow{B}_{s}\right\vert
^{2}\rightarrow 0,\text{ in probability as, }n\rightarrow \infty .%
\end{array}%
\right.
\end{equation*}%
Let the following reflected BDSDEPs with data $\left( \xi ,f,g,S\right) $%
\begin{equation*}
\left\{
\begin{array}{l}
\hat{Y}\in \mathcal{S}^{2}\left( 0,T,%
\mathbb{R}
\right) ,\text{\quad }\hat{Z}\in \mathcal{M}^{2}\left( 0,T,%
\mathbb{R}
^{d}\right) ,\text{\quad }K\in \mathcal{A}^{2},\text{\quad
}\hat{U}\in
\mathcal{L}^{2}\left( 0,T,\tilde{\mu},%
\mathbb{R}
\right) , \\
\\
\hat{Y}_{t}=\xi
+\int_{t}^{T}f(s,Y_{s},Z_{s},U_{s})ds+\int_{t}^{T}g(s,Y_{s},Z_{s},U_{s})d%
\overleftarrow{B}_{s}+\int_{t}^{T}dK_{s} \\
\\
-\int_{t}^{T}\hat{Z}_{s}dW_{s}-\int_{t}^{T}\int_{E}\hat{U}_{s}\left(
e\right) \tilde{\mu}\left( ds,de\right) , \\
\\
S_{t}\leq \hat{Y}_{t},\text{\quad }0\leq t\leq T\text{\quad and\quad }%
\int_{0}^{T}\left( \hat{Y}_{t}-S_{t}\right) dK_{t}=0.%
\end{array}%
\right.
\end{equation*}By Itô's formula, we derive that%
\begin{equation*}
\begin{array}{l}
\mathbb{E}\left\vert Y_{t}^{n}-\hat{Y}_{t}\right\vert ^{2}\leq 2\mathbb{E}%
\int_{t}^{T}\left( Y_{s}^{n}-\hat{Y}_{s}\right) \left(
f_{n}(s,Y_{s}^{n},Z_{s}^{n},U_{s}^{n})-f(s,Y_{s},Z_{s},U_{s})\right) ds \\
+2\mathbb{E}\int_{t}^{T}\left( Y_{s}^{n}-\hat{Y}_{s}\right) \left(
dK_{s}^{n}-dK_{s}\right) +\mathbb{E}\int_{t}^{T}\left\vert \left\vert
g(s,Y_{s}^{n},Z_{s}^{n},U_{s}^{n})-g(s,Y_{s},Z_{s},U_{s})\right\vert
\right\vert ^{2}ds \\
-\mathbb{E}\int_{t}^{T}\int_{E}\left\vert U_{s}^{n}\left( e\right) -\hat{U}%
_{s}\left( e\right) \right\vert ^{2}\lambda \left( de\right) ds-\mathbb{E}%
\int_{t}^{T}\left\vert Z_{s}^{n}-\hat{Z}_{s}\right\vert ^{2}ds.%
\end{array}%
\end{equation*}%
Using the fact that $\mathbb{E}\int_{t}^{T}\left( Y_{s}^{n}-\hat{Y}%
_{s}\right) \left( dK_{s}^{n}-dK_{s}\right) \leq 0,$ we get%
\begin{eqnarray*}
&&\mathbb{E}\left\vert Y_{t}^{n}-\hat{Y}_{t}\right\vert ^{2}+\mathbb{E}%
\int_{t}^{T}\int_{E}\left\vert U_{s}^{n}\left( e\right)
-\hat{U}_{s}\left( e\right) \right\vert ^{2}\lambda \left( de\right)
ds+\mathbb{E}\int_{t}^{T}\left\vert
Z_{s}^{n}-\hat{Z}_{s}\right\vert ^{2}ds \\
&\leq &2\mathbb{E}\int_{t}^{T}\left( Y_{s}^{n}-\hat{Y}_{s}\right) \left(
f_{n}(s,Y_{s}^{n},Z_{s}^{n},U_{s}^{n})-f(s,Y_{s},Z_{s},U_{s})\right) ds+%
\mathbb{E}\int_{t}^{T}\left\vert \left\vert
g(s,Y_{s}^{n},Z_{s}^{n},U_{s}^{n})-g(s,Y_{s},Z_{s},U_{s})\right\vert
\right\vert ^{2}ds,
\end{eqnarray*}%
letting $n\rightarrow \infty $, we have $Y_{t}=\hat{Y}_{t},$ $U_{t}=\hat{U}%
_{t}$ and $Z_{t}=\hat{Z}_{t}$ $d\mathbb{P}\times dt-a.e$.

Let $\left( Y^{\ast },Z^{\ast },U^{\ast },K^{\ast }\right) $ be a
solution of $\left( 1.1\right) $. Then by Theorem 3.1, we have for
any $n\in
\mathbb{N}
^{\ast }$, $Y^{n}\leq Y^{\ast }.$ Therefore, $Y$ is a minimal solution of $%
\left( 1.1\right) .$ \eop

\section{RBDSDEPs with discontinuous coefficient.}

\qquad In this section we are interested in weakening the conditions on $f$.
We assume that $f$ satisfy the following assumptions:

\vskip0.1cm\noindent \textbf{(H.8) }There exists a nonnegative process $%
f_{t}\in \mathcal{M}^{2}\left( 0,T,%
\mathbb{R}
\right) $ and constant $C>0$, such that%
\begin{equation*}
\forall \left( t,y,z,u\right) \in \left[ 0,T\right] \times \mathcal{B}%
^{2}\left(
\mathbb{R}
\right) ,\text{ }\left\vert f\left( t,y,z,u\right) \right\vert \leq
f_{t}\left( \omega \right) +C\left( \left\vert y\right\vert +\left\vert
z\right\vert +\left\vert u\right\vert \right) .
\end{equation*}%
\vskip0.1cm\noindent \textbf{(H.9) }$f\left( t,\cdot ,z,u\right) :%
\mathbb{R}
\rightarrow
\mathbb{R}
$ is a left continuous and \textbf{\ }$f\left( t,y,\cdot ,\cdot \right) $ is
a continuous.

\vskip0.1cm\noindent \textbf{(H.10) }There exists a continuous fonction%
\textbf{\ }$\pi :\left[ 0,T\right] \times \mathcal{B}^{2}\left(
\mathbb{R}
\right) $ satisfying for $y\geq y^{^{\prime }}$, $\left( z,u\right) \in
\mathbb{R}
^{d}\times L^{2}\left( E,\mathcal{E},\lambda ,%
\mathbb{R}
\right) $%
\begin{equation*}
\left\{
\begin{array}{l}
\left\vert \pi \left( t,y,z,u\right) \right\vert \leq C\left( \left\vert
y\right\vert +\left\vert z\right\vert +\left\vert u\right\vert \right) , \\
f\left( t,\omega ,y,z,u\right) -f\left( t,\omega ,y^{^{\prime }},z^{\prime
},u^{^{\prime }}\right) \geq \pi \left( t,y-y^{^{\prime }},z-z^{^{\prime
}},u-u^{^{\prime }}\right) .%
\end{array}%
\right.
\end{equation*}%
\vskip0.1cm\noindent \textbf{(H.11) }$g$ satisfies $\left(
H.2,(ii)\right) .$

\subsection{Existence result.}

The two next Lemmas will be useful in the sequel.

\begin{lemma}
\label{Lemma 5} \ Assume that $\pi $ satisfies $\left( H.10\right) ,$ $g$
satisfies $\left( H.11\right) $ and $h$ belongs in $\mathcal{M}^{2}\left( 0,T,%
\mathbb{R}
\right) $. For a continuous processes of finite variation $A$ belong in $%
\mathcal{A}^{2}$ we consider the processes $\left( \bar{Y},\bar{Z},\bar{U}%
\right) \in \mathcal{S}^{2}\left( 0,T,%
\mathbb{R}
\right) \times \mathcal{M}^{2}\left( 0,T,%
\mathbb{R}
^{d}\right) \times \mathcal{L}^{2}\left( 0,T,\tilde{\mu},%
\mathbb{R}
\right) $ such that:%
\begin{equation}
\left\{
\begin{array}{l}
\left( i\right) \text{ }\bar{Y}_{t}=\xi +\int_{t}^{T}\left( \pi \left( s,%
\bar{Y}_{s},\bar{Z}_{s},\bar{U}_{s}\right) +h\left( s\right) \right)
ds+\int_{t}^{T}g(s,\bar{Y}_{s},\bar{Z}_{s},\bar{U}_{s})d\overleftarrow{B}_{s}
\\
\\
+\int_{t}^{T}dA_{s}-\int_{t}^{T}\bar{Z}_{s}dW_{s}-\int_{t}^{T}\int_{E}\bar{U}%
_{s}\left( e\right) \tilde{\mu}\left( ds,de\right) ,\ 0\leq t\leq T, \\
\\
\left( ii\right) \text{ }\int_{0}^{T}\bar{Y}_{t}^{-}dA_{s}\geq 0.%
\end{array}%
\right.
\end{equation}%
Then we have,

$\left( 1\right) $ The \textbf{RBDSDEPs} $\left( 5.1\right) $ admits
a
minimal solution $\left( \overline{Y}_{t},\overline{Z}_{t},A_{t},\overline{U}%
_{t}\right) \in \mathcal{D}^{2}\left(
\mathbb{R}
\right) .$

$\left( 2\right) $ if $h\left( t\right) \geq 0$ and $\xi \geq 0,$ we have $%
\bar{Y}_{t}\geq 0,$ $d\mathbb{P}\times dt-a.s.$

\end{lemma}

\noindent \textbf{Proof: }$(1)$ Obtained from a previous part.

$(2)$ Applying lemma 2.1 to $\left\vert \bar{Y}_{t}^{-}\right\vert ^{2}$, we
have%
\begin{eqnarray*}
&&\mathbb{E}\left( \left\vert \bar{Y}_{t}^{-}\right\vert
^{2}+\int_{t}^{T}1_{\left\{ \bar{Y}_{s}<0\right\} }\left\vert \left\vert
\bar{Z}_{s}\right\vert \right\vert ^{2}ds+\int_{t}^{T}\int_{E}1_{\left\{
\bar{Y}_{s}<0\right\} }\left\vert \bar{U}_{s}\left( e\right) \right\vert
^{2}\lambda \left( de\right) ds\right) \\
&\leq &\mathbb{E}\left( \left\vert \xi ^{-}\right\vert ^{2}-2\int_{t}^{T}%
\bar{Y}_{s}^{-}\left( \pi \left( s,\bar{Y}_{s},\bar{Z}_{s},\bar{U}%
_{s}\right) +h\left( s\right) \right) ds-2\int_{t}^{T}\bar{Y}%
_{s}^{-}dA_{s}+\int_{t}^{T}1_{\left\{ \bar{Y}_{s}<0\right\} }\left\vert
\left\vert g(s,\bar{Y}_{s},\bar{Z}_{s},\bar{U}_{s})\right\vert \right\vert
^{2}ds\right) .
\end{eqnarray*}%
Since $h\left( t\right) \geq 0$, $\xi \geq 0$ and using the fact that $%
\int_{0}^{T}\bar{Y}_{t}^{-}dA_{s}\geq 0$, we obtain%
\begin{eqnarray*}
&&\mathbb{E}\left\vert \bar{Y}_{t}^{-}\right\vert ^{2}+\mathbb{E}%
\int_{t}^{T}1_{\left\{ \bar{Y}_{s}<0\right\} }\left\vert \left\vert \bar{Z}%
_{s}\right\vert \right\vert ^{2}ds+\mathbb{E}\int_{t}^{T}\int_{E}1_{\left\{
\bar{Y}_{s}<0\right\} }\left\vert \bar{U}_{s}\left( e\right) \right\vert
^{2}\lambda \left( de\right) ds \\
&\leq &-2\mathbb{E}\int_{t}^{T}\bar{Y}_{s}^{-}\pi \left( s,\bar{Y}_{s},\bar{Z%
}_{s},\bar{U}_{s}\right) ds+\mathbb{E}\int_{t}^{T}1_{\left\{ \bar{Y}%
_{s}<0\right\} }\left\vert \left\vert g(s,\bar{Y}_{s},\bar{Z}_{s},\bar{U}%
_{s})\right\vert \right\vert ^{2}ds.
\end{eqnarray*}%
According to assumptions $\left( H.11\right) $, we get%
\begin{eqnarray*}
&&\mathbb{E}\left\vert \bar{Y}_{t}^{-}\right\vert ^{2}+\mathbb{E}%
\int_{t}^{T}1_{\left\{ \bar{Y}_{s}<0\right\} }\left\vert \left\vert \bar{Z}%
_{s}\right\vert \right\vert ^{2}ds+\mathbb{E}\int_{t}^{T}\int_{E}1_{\left\{
\bar{Y}_{s}<0\right\} }\left\vert \bar{U}_{s}\left( e\right) \right\vert
^{2}\lambda \left( de\right) ds \\
&\leq &-2\mathbb{E}\int_{t}^{T}\bar{Y}_{s}^{-}\pi \left( s,\bar{Y}_{s},\bar{Z%
}_{s},\bar{U}_{s}\right) ds+C\mathbb{E}\int_{t}^{T}1_{\left\{ \bar{Y}%
_{s}<0\right\} }\left\vert \bar{Y}_{s}^{-}\right\vert ^{2}ds \\
&&+\alpha \mathbb{E}\int_{t}^{T}1_{\left\{ \bar{Y}_{s}<0\right\} }\left\vert
\left\vert \bar{Z}_{s}\right\vert \right\vert ^{2}ds+\alpha \mathbb{E}%
\int_{t}^{T}\int_{E}1_{\left\{ \bar{Y}_{s}<0\right\} }\left\vert \bar{U}%
_{s}\left( e\right) \right\vert ^{2}\lambda \left( de\right) ds,
\end{eqnarray*}%
applying assumption $\left( H.10\right) $ and using Young's inequality, we
have%
\begin{eqnarray*}
-2\mathbb{E}\int_{t}^{T}\bar{Y}_{s}^{-}\pi \left( s,\bar{Y}_{s},\bar{Z}_{s},%
\bar{U}_{s}\right) ds &\leq &2C\mathbb{E}\int_{t}^{T}\left\vert \bar{Y}%
_{s}^{-}\right\vert ^{2}ds+\frac{1}{2\epsilon }\mathbb{E}\int_{t}^{T}\left%
\vert \bar{Y}_{s}^{-}\right\vert ^{2}ds+2\epsilon C^{2}\mathbb{E}%
\int_{t}^{T}\left\vert \left\vert \bar{Z}_{s}\right\vert \right\vert ^{2}ds
\\
&&+\frac{1}{2\epsilon }\mathbb{E}\int_{t}^{T}\left\vert \bar{Y}%
_{s}^{-}\right\vert ^{2}ds+2\epsilon C^{2}\mathbb{E}\int_{t}^{T}\int_{E}%
\left\vert \bar{U}_{s}\left( e\right) \right\vert ^{2}\lambda \left(
de\right) ds.
\end{eqnarray*}%
Then%
\begin{eqnarray*}
&&\mathbb{E}\left\vert \bar{Y}_{t}^{-}\right\vert ^{2}+\mathbb{E}%
\int_{t}^{T}1_{\left\{ \bar{Y}_{s}<0\right\} }\left\vert \left\vert \bar{Z}%
_{s}\right\vert \right\vert ^{2}ds+\mathbb{E}\int_{t}^{T}\int_{E}1_{\left\{
\bar{Y}_{s}<0\right\} }\left\vert \bar{U}_{s}\left( e\right) \right\vert
^{2}\lambda \left( de\right) ds \\
&\leq &\left( 3C+\epsilon ^{-1}\right) \mathbb{E}\int_{t}^{T}\left\vert \bar{%
Y}_{s}^{-}\right\vert ^{2}ds+\left( \alpha +2\epsilon C^{2}\right) \mathbb{E}%
\int_{t}^{T}1_{\left\{ \bar{Y}_{s}<0\right\} }\left( \left\vert \left\vert
\bar{Z}_{s}\right\vert \right\vert ^{2}+\int_{E}\left\vert \bar{U}_{s}\left(
e\right) \right\vert ^{2}\lambda \left( de\right) \right) ds
\end{eqnarray*}%
Therefore, choosing $\epsilon $, $\alpha $ and $C$ such that $0<\alpha
+2\epsilon C^{2}<1$ and using Gronwall's inequality, we have%
\begin{equation*}
\mathbb{E}\left\vert \bar{Y}_{t}^{-}\right\vert ^{2}=0,
\end{equation*}%
$\mathbf{P}-a.s.$ for all $t\in \left[ 0,T\right] .$ Finally implies that $%
\bar{Y}_{t}\geq 0,$ $\mathbf{P}-a.s.$ for all $t\in \left[ 0,T\right] .$\eop

Now by theorem (4.1), we consider the processes $\left( \tilde{Y}_{t}^{0},\tilde{Z}%
_{t}^{0},\tilde{K}_{t}^{0},\tilde{U}_{t}^{0}\right) ,$ $\left(
Y_{t}^{0},Z_{t}^{0},K_{t}^{0},U_{t}^{0}\right) $ and sequence of processes $%
\left( \tilde{Y}_{t}^{n},\tilde{Z}_{t}^{n},\tilde{K}_{t}^{n},\tilde{U}%
_{t}^{n}\right) _{n\geq 0}$ respectively minimal solution of the following
RBDSDEPs for all $t\in \left[ 0,T\right] $%
\begin{equation}
\left\{
\begin{array}{l}
\left( i\right) \tilde{Y}_{t}^{0}=\xi +\int_{t}^{T}\left[ -C\left(
\left\vert \tilde{Y}_{s}^{0}\right\vert +\left\vert \tilde{Z}%
_{s}^{0}\right\vert +\left\vert \tilde{U}_{s}^{0}\right\vert \right) -f_{s}%
\right] ds+\int_{t}^{T}g(s,\tilde{Y}_{s}^{0},\tilde{Z}_{s}^{0},\tilde{U}%
_{s}^{0})d\overleftarrow{B}_{s} \\
\\
+\int_{t}^{T}d\tilde{K}_{s}^{0}-\int_{t}^{T}\tilde{Z}_{s}^{0}dW_{s}-%
\int_{t}^{T}\int_{E}\tilde{U}_{s}^{0}\left( e\right) \tilde{\mu}\left(
ds,de\right) ,\ 0\leq t\leq T, \\
\\
\left( ii\right)\tilde{Y}_{t}^{0}\geq S_{t}, \\
\\
\left( iii\right) \int_{0}^{T}\left( \tilde{Y}_{s}^{0}-S_{s}\right) d%
\tilde{K}_{s}^{0}=0,%
\end{array}%
\right.
\end{equation}%
\begin{equation}
\left\{
\begin{array}{l}
\left( i\right) \text{ }Y_{t}^{0}=\xi +\int_{t}^{T}\left[ C\left( \left\vert
Y_{s}^{0}\right\vert +\left\vert Z_{s}^{0}\right\vert +\left\vert
U_{s}^{0}\right\vert \right) +f_{s}\right] ds+%
\int_{t}^{T}g(s,Y_{s}^{0},Z_{s}^{0},U_{s}^{0})d\overleftarrow{B}_{s} \\
\\
+\int_{t}^{T}dK_{s}^{0}-\int_{t}^{T}Z_{s}^{0}dW_{s}-\int_{t}^{T}%
\int_{E}U_{s}^{0}\left( e\right) \tilde{\mu}\left( ds,de\right) ,\ 0\leq
t\leq T, \\
\\
\left( ii\right) \text{ }Y_{t}^{0}\geq S_{t}, \\
\\
\left( iii\right) \text{ }\int_{0}^{T}\left( Y_{s}^{0}-S_{s}\right)
dK_{s}^{0}=0,%
\end{array}%
\right.
\end{equation}
and
\begin{equation}
\left\{
\begin{array}{l}
\left( i\right) \text{ }\tilde{Y}_{t}^{n}=\xi +\int_{t}^{T}\left[ f(s,\tilde{%
Y}_{s}^{n-1},\tilde{Z}_{s}^{n-1},\tilde{U}_{s}^{n-1})ds+\pi \left( s,\tilde{Y%
}_{t}^{n}-\tilde{Y}_{t}^{n-1},\tilde{Z}_{t}^{n}-\tilde{Z}_{t}^{n-1},\tilde{U}%
_{t}^{n}-\tilde{U}_{t}^{n-1}\right) \right] ds \\
\\
+\int_{t}^{T}g(s,\tilde{Y}_{s}^{n},\tilde{Z}_{s}^{n},\tilde{U}_{s}^{n})d%
\overleftarrow{B}_{s}+\int_{t}^{T}d\tilde{K}_{s}^{n}-\int_{t}^{T}\tilde{Z}%
_{s}^{n}dW_{s}-\int_{t}^{T}\int_{E}\tilde{U}_{s}^{n}\left( e\right) \tilde{%
\mu}\left( ds,de\right) ,\ 0\leq t\leq T, \\
\\
\left( ii\right) \tilde{Y}_{t}^{n}\geq S_{t}, \\
\\
\left( iii\right)\int_{0}^{T}\left( \tilde{Y}_{s}^{n}-S_{s}\right) d%
\tilde{K}_{s}^{n}=0.%
\end{array}%
\right.
\end{equation}

\begin{lemma}
\label{Lemma 6}Under the assumptions $( H.1),\ ( H.3),\ ( H.5) $
and $\left( H.8\right) -\left( H.11\right) $, we have for any $n\geq 1$ and $%
t\in \left[ 0,T\right] $%
\begin{equation*}
\tilde{Y}_{t}^{0}\leq \tilde{Y}_{t}^{n}\leq \tilde{Y}_{t}^{n+1}\leq
Y_{t}^{0}.
\end{equation*}%
%
%
%
%
%
%
%
%
\end{lemma}

\noindent \textbf{Proof: }For any $n\geq 0$, we set%
\begin{equation*}
\left\{
\begin{array}{l}
\delta \rho _{t}^{n+1}=\rho _{t}^{n+1}-\rho _{t}^{n}, \\
and \\
\Delta \psi ^{n+1}(s,\delta \tilde{Y}_{s}^{n+1},\delta \tilde{Z}%
_{s}^{n+1},\delta \tilde{U}_{s}^{n+1})=\psi (s,\delta \tilde{Y}_{s}^{n+1}+%
\tilde{Y}_{s}^{n},\delta \tilde{Z}_{s}^{n+1}+\tilde{Z}_{s}^{n},\delta \tilde{%
U}_{s}^{n+1}+\tilde{U}_{s}^{n})-\psi (s,\tilde{Y}_{s}^{n},\tilde{Z}_{s}^{n},%
\tilde{U}_{s}^{n}).%
\end{array}%
\right.
\end{equation*}%
Using Eq.$\left( 5.4\right) $, we have%
\begin{eqnarray*}
\delta \tilde{Y}_{t}^{n+1} &=&\int_{t}^{T}\left[ \pi \left( s,\delta \tilde{Y%
}_{s}^{n+1},\delta \tilde{Z}_{s}^{n+1},\delta \tilde{U}_{s}^{n+1}\right)
+\theta _{s}^{n}\right] ds+\int_{t}^{T}\Delta g^{n+1}(s,\delta \tilde{Y}%
_{s}^{n+1},\delta \tilde{Z}_{s}^{n+1},\delta \tilde{U}_{s}^{n+1})d%
\overleftarrow{B}_{s} \\
&&+\int_{t}^{T}d\left( \delta \tilde{K}_{s}^{n+1}\right)
-\int_{t}^{T}\delta
\tilde{Z}_{s}^{n+1}dW_{s}-\int_{t}^{T}\int_{E}\delta \tilde{U}%
_{s}^{n+1}\left( e\right) \tilde{\mu}\left( ds,de\right) ,
\end{eqnarray*}%
where%
\begin{equation*}
\left\{
\begin{array}{l}
\theta _{s}^{n}=\Delta f^{n}(s,\delta\tilde{Y}_{s}^{n},\delta\tilde{Z}_{s}^{n},%
\delta\tilde{U}_{s}^{n})-\pi \left( s,\delta \tilde{Y}_{s}^{n},\delta \tilde{Z}%
_{s}^{n},\delta \tilde{U}_{s}^{n}\right) , \\
and \\
\theta _{s}^{0}=f(s,\tilde{Y}_{s}^{0},\tilde{Z}_{s}^{0},\tilde{U}%
_{s}^{0})+C\left( \left\vert \tilde{Y}_{s}^{0}\right\vert +\left\vert \tilde{%
Z}_{s}^{0}\right\vert +\left\vert \tilde{U}_{s}^{0}\right\vert \right)
+f_{s},\forall n\geq 0.%
\end{array}%
\right.
\end{equation*}%
According to it definition, on can show that $\theta _{s}^{0}$ and
$\Delta g^{n}$, $\forall n\geq 0$ satisfy all assumption of lemma
5.1. Moreover, since $\tilde{K}_{t}^{n}$ is a continuous and
increasing process, for all $ n\geq 0$, $\delta \tilde{K}_{s}^{n+1}$
is a continuous process of finite variation. Using the same argument
as in first part, On can show that
\begin{eqnarray*}
\int_{0}^{T}\left( \tilde{Y}_{t}^{n+1}-\tilde{Y}_{t}^{n}\right)^{-}
d\left( \delta \tilde{K}_{t}^{n+1}\right) &=&\int_{0}^{T}\left(
\tilde{Y}_{t}^{n+1}-\tilde{Y}_{t}^{n}\right)^{-} d\tilde{K}
_{t}^{n+1}-\int_{0}^{T}\left(
\tilde{Y}_{t}^{n+1}-\tilde{Y}_{t}^{n}\right)^{-} d
\tilde{K}_{t}^{n}\geq 0.
\end{eqnarray*}%
Applying lemma 5.1 we deduce that $\delta \tilde{Y}_{t}^{n+1}\geq 0,$ i.e. $%
\tilde{Y}_{t}^{n+1}\geq \tilde{Y}_{t}^{n}$ , for all $t\in \left[
0,T\right] ,$ we have $\tilde{Y}_{t}^{n+1}\geq \tilde{Y}_{t}^{n}\geq
\tilde{Y}_{t}^{0}.$

Now we show prove that $\tilde{Y}_{t}^{n+1}\leq Y_{t}^{0},$ by
definition, we obtain

\begin{eqnarray*}
Y_{t}^{0}-\tilde{Y}_{t}^{n+1}
&=&\int_{t}^{T}\left( -C\left( \left\vert Y_{s}^{0}-\tilde{Y}%
_{s}^{n+1}\right\vert +\left\vert Z_{s}^{0}-\tilde{Z}_{s}^{n+1}\right\vert
+\left\vert U_{s}^{0}-\tilde{U}_{s}^{n+1}\right\vert \right) +\Lambda
_{s}^{n+1}\right) ds \\
&&+\int_{t}^{T}\left( g(s,Y_{s}^{0},Z_{s}^{0},U_{s}^{0})-g(s,\tilde{Y}%
_{s}^{n+1},\tilde{Z}_{s}^{n+1},\tilde{U}_{s}^{n+1})\right) d\overleftarrow{B}%
_{s} \\
&&+\int_{t}^{T}\left( dK_{s}^{0}-d\tilde{K}_{s}^{n+1}\right)
+\int_{t}^{T}\left( Z_{s}^{0}-\tilde{Z}_{s}^{n+1}\right)
dW_{s}-\int_{t}^{T}\int_{E}\left( U_{s}^{0}\left( e\right) -\tilde{U}%
_{s}^{n+1}\left( e\right) \right) \tilde{\mu}\left( ds,de\right) ,
\end{eqnarray*}%
where

$\Lambda _{s}^{n+1}=C\left( \left\vert Y_{s}^{0}-\tilde{Y}%
_{s}^{n+1}\right\vert +\left\vert Z_{s}^{0}-\tilde{Z}_{s}^{n+1}\right\vert
+\left\vert U_{s}^{0}-\tilde{U}_{s}^{n+1}\right\vert +\left\vert
Y_{s}^{0}\right\vert +\left\vert Z_{s}^{0}\right\vert +\left\vert
U_{s}^{0}\right\vert \right) +f_{s}-f(s,\tilde{Y}_{s}^{n},\tilde{Z}_{s}^{n},%
\tilde{U}_{s}^{n})-\pi \left( s,\delta \tilde{Y}_{s}^{n+1},\delta \tilde{Z}%
_{s}^{n+1},\delta \tilde{U}_{s}^{n+1}\right) .$

Also using lemma 5.1 we deduce that
$Y_{t}^{0}-\tilde{Y}_{t}^{n+1}\geq 0$, i.e. $Y_{t}^{0}\geq
\tilde{Y}_{t}^{n+1},$ for all $t\in \left[ 0,T\right] $.
Thus, we have for all $n\geq 0$%
\begin{equation*}
Y_{t}^{0}\geq \tilde{Y}_{t}^{n+1}\geq \tilde{Y}_{t}^{n}\geq \tilde{Y}%
_{t}^{0},\text{ }d\mathbb{P}\times dt-a.s.,\text{\quad }\forall t\in \left[
0,T\right] .
\end{equation*}%
\eop Now, our main result
\begin{theorem}
\label{th1 copy(1)}Under assumption $\left( H.1\right),\ \left(
H.3\right),\ \left( H.5\right) $ and $\left( H.8\right) -\left(
H.11\right) $, the RBDSDEPs $\left( 1.1\right) $ has a minimal
solution $(Y_{t},Z_{t},K_{t},U_{t})_{0\leq t\leq T}\in
\mathcal{D}^{2}\left(
\mathbb{R}
\right) .$

\end{theorem}

\noindent \textbf{Proof: }Since $\left\vert \tilde{Y}_{t}^{n}\right\vert
\leq \max \left( \tilde{Y}_{t}^{0},Y_{t}^{0}\right) \leq \left\vert \tilde{Y}%
_{t}^{0}\right\vert +\left\vert Y_{t}^{0}\right\vert $ for all $t\in \left[
0,T\right] $, we have%
\begin{equation*}
\sup_{n}\mathbb{E}\left( \sup_{0\leq t\leq T}\left\vert \tilde{Y}%
_{t}^{n}\right\vert ^{2}\right) \leq \mathbb{E}\left( \sup_{0\leq t\leq
T}\left\vert \tilde{Y}_{t}^{0}\right\vert ^{2}\right) +\mathbb{E}\left(
\sup_{0\leq t\leq T}\left\vert Y_{t}^{0}\right\vert ^{2}\right) <\infty .
\end{equation*}%
Therefore, we deduce from the Lebesgue's dominated convergence theorem that $%
\left( \tilde{Y}_{t}^{n}\right) _{n\geq 0}$ converges in $\mathcal{S}%
^{2}\left( 0,T,%
\mathbb{R}
\right) $ to a limit $Y.$

On the other hand from $(5.4)$, we deduce that%
\begin{eqnarray*}
\tilde{Y}_{0}^{n+1} &=&\tilde{Y}_{T}^{n+1}+\int_{0}^{T}\left[ f(s,\tilde{Y}%
_{s}^{n},\tilde{Z}_{s}^{n},\tilde{U}_{s}^{n})ds+\pi \left( s,\delta \tilde{Y}%
_{t}^{n+1},\delta \tilde{Z}_{t}^{n+1},\delta \tilde{U}_{t}^{n+1}\right) %
\right] ds \\
&&+\int_{0}^{T}g(s,\tilde{Y}_{s}^{n+1},\tilde{Z}_{s}^{n+1},\tilde{U}%
_{s}^{n+1})d\overleftarrow{B}_{s}+\int_{0}^{T}d\tilde{K}_{s}^{n+1}-%
\int_{0}^{T}\tilde{Z}_{s}^{n+1}dW_{s}-\int_{0}^{T}\int_{E}\tilde{U}%
_{s}^{n+1}\left( e\right) \tilde{\mu}\left( ds,de\right) ,
\end{eqnarray*}%
applying lemma 2.1, we obtain%
\begin{eqnarray*}
&&\mathbb{E}\left\vert \tilde{Y}_{0}^{n+1}\right\vert ^{2}+\mathbb{E}%
\int_{0}^{T}\left\vert \tilde{Z}_{s}^{n+1}\right\vert ^{2}ds+\mathbb{E}%
\int_{0}^{T}\int_{E}\left\vert \tilde{U}_{s}^{n+1}\left( e\right)
\right\vert ^{2}\lambda \left( de\right) ds \\
&\leq &\mathbb{E}\left\vert \xi \right\vert ^{2}+2\mathbb{E}\int_{0}^{T}%
\tilde{Y}_{s}^{n+1}\left( f(s,\tilde{Y}_{s}^{n},\tilde{Z}_{s}^{n},\tilde{U}%
_{s}^{n})+\pi \left( s,\delta \tilde{Y}_{s}^{n+1},\delta \tilde{Z}%
_{s}^{n+1},\delta \tilde{U}_{s}^{n+1}\right) \right) ds \\
&&+2\mathbb{E}\int_{0}^{T}\tilde{Y}_{s}^{n+1}d\tilde{K}_{s}^{n+1}+\mathbb{E}\int_{0}^{T}\left%
\vert \left\vert g(s,\tilde{Y}_{s}^{n+1},\tilde{Z}_{s}^{n+1},\tilde{U}%
_{s}^{n+1})\right\vert \right\vert ^{2}ds.
\end{eqnarray*}%
From $\left( H.8\right) $ and $\left( H.10\right) $, we get%
\begin{eqnarray*}
&&\tilde{Y}_{t}^{n+1}\left( f(t,\tilde{Y}_{t}^{n},\tilde{Z}_{t}^{n},\tilde{U}%
_{t}^{n})+\pi \left( t,\delta \tilde{Y}_{t}^{n+1},\delta \tilde{Z}%
_{t}^{n+1},\delta \tilde{U}_{t}^{n+1}\right) \right)  \\
&\leq &\left\vert \tilde{Y}_{t}^{n+1}\right\vert \left\{ f_{t}\left( \omega
\right) +2C\left( \left\vert \tilde{Y}_{t}^{n}\right\vert +\left\vert \tilde{%
Z}_{t}^{n}\right\vert +\left\vert \tilde{U}_{t}^{n}\right\vert \right)
+C\left( \left\vert \tilde{Y}_{t}^{n+1}\right\vert +\left\vert \tilde{Z}%
_{t}^{n+1}\right\vert +\left\vert \tilde{U}_{t}^{n+1}\right\vert \right)
\right\}  \\
&\leq &\frac{\left\vert \tilde{Y}_{t}^{n+1}\right\vert ^{2}}{2}+\frac{%
f_{t}\left( \omega \right) }{2}+C^{2}\left\vert \tilde{Y}_{t}^{n+1}\right%
\vert ^{2}+\left\vert \tilde{Y}_{t}^{n}\right\vert ^{2}+\frac{2C^{2}}{%
\epsilon _{1}}\left\vert \tilde{Y}_{t}^{n+1}\right\vert ^{2}+\frac{\epsilon
_{1}}{2}\left\vert \tilde{Z}_{t}^{n}\right\vert ^{2}+\frac{2C^{2}}{\epsilon
_{2}}\left\vert \tilde{Y}_{t}^{n+1}\right\vert ^{2}+\frac{\epsilon _{2}}{2}%
\left\vert \tilde{U}_{t}^{n}\right\vert ^{2} \\
&&+C\left\vert \tilde{Y}_{t}^{n+1}\right\vert ^{2}+\frac{C^{2}}{2\epsilon
_{3}}\left\vert \tilde{Y}_{t}^{n+1}\right\vert ^{2}+\frac{\epsilon _{3}}{2}%
\left\vert \tilde{Z}_{t}^{n+1}\right\vert ^{2}+\frac{C^{2}}{2\epsilon _{4}}%
\left\vert \tilde{Y}_{t}^{n+1}\right\vert ^{2}+\frac{\epsilon _{4}}{2}%
\left\vert \tilde{U}_{t}^{n+1}\right\vert ^{2}. \\
&=&\left( \frac{1}{2}+C^{2}+\frac{2C^{2}}{\epsilon _{1}}+\frac{2C^{2}}{%
\epsilon _{2}}+\frac{C^{2}}{2\epsilon _{3}}+\frac{C^{2}}{2\epsilon _{4}}%
+C\right) \left\vert \tilde{Y}_{t}^{n+1}\right\vert ^{2} \\
&&+\frac{\epsilon _{3}}{2}\left\vert \tilde{Z}_{t}^{n+1}\right\vert ^{2}+%
\frac{\epsilon _{4}}{2}\left\vert \tilde{U}_{t}^{n+1}\right\vert
^{2}+\left\vert \tilde{Y}_{t}^{n}\right\vert ^{2}+\frac{\epsilon _{1}}{2}%
\left\vert \tilde{Z}_{t}^{n}\right\vert ^{2}+\frac{\epsilon _{2}}{2}%
\left\vert \tilde{U}_{t}^{n}\right\vert ^{2}+\frac{f_{t}\left( \omega
\right) }{2}.
\end{eqnarray*}%
Also applying $\left( H.11\right) $, we obtain the following inequality%
\begin{eqnarray*}
\left\vert \left\vert g(s,\tilde{Y}_{s}^{n+1},\tilde{Z}_{s}^{n+1},\tilde{U}%
_{s}^{n+1})\right\vert \right\vert ^{2} &\leq &\left\vert \left\vert g(s,%
\tilde{Y}_{s}^{n+1},\tilde{Z}_{s}^{n+1},\tilde{U}_{s}^{n+1})-g(s,0,0,0)%
\right\vert \right\vert ^{2}+\left\vert \left\vert g(s,0,0,0)\right\vert
\right\vert ^{2}, \\
&\leq &C\left\vert \tilde{Y}_{s}^{n+1}\right\vert ^{2}+\alpha \left\{
\left\vert \tilde{Z}_{s}^{n+1}\right\vert ^{2}+\left\vert \tilde{U}%
_{s}^{n+1}\right\vert ^{2}\right\} +\left\vert \left\vert
g(s,0,0,0)\right\vert \right\vert ^{2}.
\end{eqnarray*}%
Using Young's inequality, we get%
\begin{equation*}
2\mathbb{E}\int_{0}^{T}\tilde{Y}_{s}^{n+1}d\tilde{K}_{s}^{n+1}\leq 2\mathbb{E%
}\int_{0}^{T}S_{s}d\tilde{K}_{s}^{n+1}\leq \frac{1}{\theta
}\mathbb{E}\left(
\sup_{0\leq t\leq T}\left\vert S_{t}\right\vert ^{2}\right) +\theta \mathbb{E%
}\left\vert \tilde{K}_{T}^{n+1}\right\vert ^{2}.
\end{equation*}%
Therefore, there exists a constant $C$ independent of $n$ such that for any $%
\epsilon _{i},$ where $i=1:4$, we derive%
\begin{eqnarray}
&&\mathbb{E}\int_{0}^{T}\left\vert \tilde{Z}_{s}^{n+1}\right\vert ^{2}ds+%
\mathbb{E}\int_{0}^{T}\int_{E}\left\vert \tilde{U}_{s}^{n+1}\left( e\right)
\right\vert ^{2}\lambda \left( de\right) ds  \notag \\
&\leq &C+\left( \epsilon _{3}+\alpha \right) \mathbb{E}\int_{0}^{T}\left%
\vert \tilde{Z}_{s}^{n+1}\right\vert ^{2}ds+\left( \epsilon _{4}+\alpha
\right) \mathbb{E}\int_{0}^{T}\int_{E}\left\vert \tilde{U}%
_{s}^{n+1}\right\vert ^{2}\lambda \left( de\right) ds \\
&&+\epsilon _{1}\mathbb{E}\int_{0}^{T}\left\vert \tilde{Z}%
_{s}^{n}\right\vert ^{2}ds+\epsilon _{2}\mathbb{E}\int_{t}^{T}\int_{E}\left%
\vert \tilde{U}_{s}^{n}\right\vert ^{2}\lambda \left( de\right)
ds+\theta \mathbb{E}\left\vert \tilde{K}_{T}^{n+1}\right\vert ^{2}.
\notag
\end{eqnarray}%
Moreover, since%
\begin{eqnarray*}
\tilde{K}_{T}^{n+1} &=&\tilde{Y}_{0}^{n+1}-\xi -\int_{0}^{T}\left[ f(s,%
\tilde{Y}_{s}^{n},\tilde{Z}_{s}^{n},\tilde{U}_{s}^{n})ds+\pi \left( s,\delta
\tilde{Y}_{s}^{n+1},\delta \tilde{Z}_{s}^{n+1},\delta \tilde{U}%
_{s}^{n+1}\right) \right] ds \\
&&-\int_{0}^{T}g(s,\tilde{Y}_{s}^{n+1},\tilde{Z}_{s}^{n+1},\tilde{U}%
_{s}^{n+1})d\overleftarrow{B}_{s}+\int_{0}^{T}\tilde{Z}_{s}^{n+1}dW_{s}+%
\int_{0}^{T}\int_{E}\tilde{U}_{s}^{n+1}\left( e\right) \tilde{\mu}\left(
ds,de\right) ,
\end{eqnarray*}%
Using Hölder's inequality and assumption $\left( H.8\right) ,$ $\left(
H.10\right) $, there exists two constants $C_1$ and $C_2$ depend of $\xi, C,\alpha, \epsilon_i, i=1,...,4$, and we have that%
\begin{equation*}
\mathbb{E}\left\vert \tilde{K}_{T}^{n+1}\right\vert ^{2}\leq
C_{1}+C_{2}\left( \mathbb{E}\int_{0}^{T}\left( \left\vert \tilde{Z}%
_{s}^{n}\right\vert ^{2}+\left\vert \tilde{Z}_{s}^{n+1}\right\vert
^{2}\right) ds+\mathbb{E}\int_{0}^{T}\int_{E}\left( \left\vert \tilde{U}%
_{s}^{n}\right\vert ^{2}+\left\vert \tilde{U}_{s}^{n+1}\right\vert
^{2}\right) \lambda \left( de\right) ds\right) ,
\end{equation*}%
we come back to inequality $\left( 5.5\right) $, we obtain%
\begin{eqnarray*}
&&\mathbb{E}\int_{0}^{T}\left\vert \tilde{Z}_{s}^{n+1}\right\vert ^{2}ds+%
\mathbb{E}\int_{0}^{T}\int_{E}\left\vert \tilde{U}_{s}^{n+1}\left( e\right)
\right\vert ^{2}\lambda \left( de\right) ds \\
&\leq &\left( C+\theta C_{1}\right) +\left( \epsilon _{1}+\theta
C_{2}\right) \mathbb{E}\int_{0}^{T}\left\vert \tilde{Z}_{s}^{n}\right\vert
^{2}ds+\left( \epsilon _{2}+\theta C_{2}\right) \mathbb{E}%
\int_{0}^{T}\int_{E}\left\vert \tilde{U}_{s}^{n}\right\vert ^{2}\lambda
\left( de\right) ds \\
&&+\left( \epsilon _{3}+\alpha +\theta C_{2}\right) \mathbb{E}%
\int_{0}^{T}\left\vert \tilde{Z}_{s}^{n+1}\right\vert ^{2}ds+\left( \epsilon
_{4}+\alpha +\theta C_{2}\right) \mathbb{E}\int_{0}^{T}\int_{E}\left\vert
\tilde{U}_{s}^{n+1}\right\vert ^{2}\lambda \left( de\right) ds,
\end{eqnarray*}%
we taking $\epsilon _{1}=\epsilon _{2}=\epsilon _{0}$ and $\epsilon
_{3}=\epsilon _{4}=\bar{\epsilon},$ we have%
\begin{eqnarray*}
&&\mathbb{E}\int_{0}^{T}\left\vert \tilde{Z}_{s}^{n+1}\right\vert ^{2}ds+%
\mathbb{E}\int_{0}^{T}\int_{E}\left\vert \tilde{U}_{s}^{n+1}\left( e\right)
\right\vert ^{2}\lambda \left( de\right) ds \\
&\leq &\left( C+\theta C_{1}\right) +\left( \epsilon _{0}+\theta
C_{2}\right) \left\{ \mathbb{E}\int_{0}^{T}\left\vert \tilde{Z}%
_{s}^{n}\right\vert ^{2}ds+\mathbb{E}\int_{0}^{T}\int_{E}\left\vert \tilde{U}%
_{s}^{n}\right\vert ^{2}\lambda \left( de\right) ds\right\}  \\
&&+\left( \bar{\epsilon}+\theta C_{2}+\alpha \right) \mathbb{E}%
\int_{0}^{T}\left( \left\vert \tilde{Z}_{s}^{n+1}\right\vert
^{2}+\int_{E}\left\vert \tilde{U}_{s}^{n+1}\left( e\right) \right\vert
^{2}\lambda \left( de\right) \right) ds,
\end{eqnarray*}%
we chossing $\bar{\epsilon},$ $\theta $ and $\alpha $ such that $0\leq
\left( \bar{\epsilon}+\theta C_{2}+\alpha \right) <1,$ we get%
\begin{eqnarray*}
&&\mathbb{E}\int_{0}^{T}\left\vert \tilde{Z}_{s}^{n+1}\right\vert ^{2}ds+%
\mathbb{E}\int_{0}^{T}\int_{E}\left\vert \tilde{U}_{s}^{n+1}\left( e\right)
\right\vert ^{2}\lambda \left( de\right) ds \\
&\leq &\left( C+\theta C_{1}\right) +\left( \epsilon _{0}+\theta
C_{2}\right) \left\{ \mathbb{E}\int_{0}^{T}\left\vert \tilde{Z}%
_{s}^{n}\right\vert ^{2}ds+\mathbb{E}\int_{0}^{T}\int_{E}\left\vert \tilde{U}%
_{s}^{n}\right\vert ^{2}\lambda \left( de\right) ds\right\}  \\
&\leq &\left( C+\theta C_{1}\right) \sum_{i=0}^{i=n-1}\left( \epsilon
_{0}+\theta C_{2}\right) ^{i}+\left( \epsilon _{0}+\theta C_{2}\right)
^{n}\left\{ \mathbb{E}\int_{0}^{T}\left\vert \tilde{Z}_{s}^{0}\right\vert
^{2}ds+\mathbb{E}\int_{0}^{T}\int_{E}\left\vert \tilde{U}_{s}^{0}\right\vert
^{2}\lambda \left( de\right) ds\right\} .,
\end{eqnarray*}%
Now chossing $\epsilon _{0},$ $\theta $ and $C_{2}$ such that $\epsilon
_{0}+\theta C_{2}<1$ and notting $\mathbb{E}\int_{0}^{T}\left( \left\vert
\tilde{Z}_{s}^{0}\right\vert ^{2}+\int_{E}\left\vert \tilde{U}%
_{s}^{0}\right\vert ^{2}\lambda \left( de\right) \right) ds<\infty .$ We obtain%
\begin{equation*}
\begin{tabular}{lll}
$\sup_{n\in
\mathbb{N}
}\mathbb{E}\int_{0}^{T}\left\vert \tilde{Z}_{s}^{n+1}\right\vert
^{2}ds<\infty $ & and & $\sup_{n\in
\mathbb{N}
}\mathbb{E}\int_{0}^{T}\int_{E}\left\vert \tilde{U}_{s}^{n+1}\left( e\right)
\right\vert ^{2}\lambda \left( de\right) ds<\infty ,$%
\end{tabular}%
\end{equation*}%
consequently, we deduce that%
\begin{equation*}
\mathbb{E}\left\vert \tilde{K}_{T}^{n+1}\right\vert ^{2}<\infty .
\end{equation*}%
Now we shall prove that $\left( \tilde{Z}^{n},\tilde{K}^{n},\tilde{U}%
^{n}\right) $ is a Cauchy sequence in $\mathcal{M}^{2}\left( 0,T,%
\mathbb{R}
^{d}\right) \times \mathcal{A}^{2}\times \mathcal{L}^{2}\left( 0,T,\tilde{\mu%
},%
\mathbb{R}
\right) ,$ set $\Gamma _{s}^{n}=f(s,\tilde{Y}_{s}^{n-1},\tilde{Z}_{s}^{n-1},%
\tilde{U}_{s}^{n-1})+\pi \left( s,\delta \tilde{Y}_{s}^{n},\delta \tilde{Z}%
_{s}^{n},\delta \tilde{U}_{s}^{n}\right) ,$ we have

\begin{eqnarray*}
\tilde{Y}_{t}^{n}-\tilde{Y}_{t}^{m} &=&\int_{t}^{T}\left( \Gamma
_{s}^{n}-\Gamma _{s}^{m}\right) ds+\int_{t}^{T}\left( g(s,\tilde{Y}_{s}^{n},%
\tilde{Z}_{s}^{n},\tilde{U}_{s}^{n})-g(s,\tilde{Y}_{s}^{m},\tilde{Z}_{s}^{m},%
\tilde{U}_{s}^{m})\right) d\overleftarrow{B}_{s} \\
&&+\int_{t}^{T}\left( d\tilde{K}_{s}^{n}-d\tilde{K}_{s}^{m}\right)
-\int_{t}^{T}\left( \tilde{Z}_{s}^{n}-\tilde{Z}_{s}^{m}\right)
dW_{s}-\int_{t}^{T}\int_{E}\left( \tilde{U}_{s}^{n}\left( e\right) -\tilde{U}%
_{s}^{m}\left( e\right) \right) \tilde{\mu}\left( ds,de\right) ,
\end{eqnarray*}%
applying Lemma 2.1 to $\left\vert \delta \tilde{Y}_{s}^{n,m}\right\vert
^{2}=\left\vert \tilde{Y}_{s}^{n}-\tilde{Y}_{s}^{m}\right\vert ^{2}$, we have%
\begin{eqnarray*}
&&\mathbb{E}\left\vert \tilde{Y}_{t}^{n}-\tilde{Y}_{t}^{m}\right\vert ^{2}+%
\mathbb{E}\int_{t}^{T}\left\vert \tilde{Z}_{s}^{n}-\tilde{Z}%
_{s}^{m}\right\vert ^{2}ds+\mathbb{E}\int_{t}^{T}\int_{E}\left\vert \tilde{U}%
_{s}^{n}-\tilde{U}_{s}^{m}\right\vert ^{2}\lambda \left( de\right) ds \\
&\leq &2\mathbb{E}\int_{t}^{T}\left( \tilde{Y}_{s}^{n}-\tilde{Y}
_{s}^{m}\right) \left( \Gamma _{s}^{n}-\Gamma _{s}^{m}\right)
ds+2\mathbb{E}\int_{t}^{T}\left( \tilde{Y}_{s}^{n+1}-\tilde{Y}_{s}^{n}\right) \left( d%
\tilde{K}_{s}^{n}-d\tilde{K}_{s}^{m}\right)  \\
&&+\mathbb{E}\int_{t}^{T}\left\vert \left\vert \left( g(s,\tilde{Y}_{s}^{n},\tilde{Z}%
_{s}^{n},\tilde{U}_{s}^{n})-g(s,\tilde{Y}_{s}^{m},\tilde{Z}_{s}^{m},\tilde{U}%
_{s}^{m})\right) \right\vert \right\vert ^{2}ds,
\end{eqnarray*}%
since $\int_{t}^{T}\left( \tilde{Y}_{s}^{n+1}-\tilde{Y}_{s}^{n}\right)
\left( d\tilde{K}_{s}^{n}-d\tilde{K}_{s}^{m}\right) \leq 0,$ we obtain%
\begin{eqnarray*}
&&\mathbb{E}\int_{0}^{T}\left\vert \tilde{Z}_{s}^{n}-\tilde{Z}%
_{s}^{m}\right\vert ^{2}ds+\mathbb{E}\int_{t}^{T}\int_{E}\left\vert \tilde{U}%
_{s}^{n}-\tilde{U}_{s}^{m}\right\vert ^{2}\lambda \left( de\right) ds \\
&\leq &2\mathbb{E}\int_{t}^{T}\left( \tilde{Y}_{s}^{n}-\tilde{Y}%
_{s}^{m}\right) \left( \Gamma _{s}^{n}-\Gamma _{s}^{m}\right)
ds+\mathbb{E}\int_{t}^{T}\left\vert \left\vert \left( g(s,\tilde{Y}_{s}^{n},\tilde{Z}%
_{s}^{n},\tilde{U}_{s}^{n})-g(s,\tilde{Y}_{s}^{m},\tilde{Z}_{s}^{m},\tilde{U}%
_{s}^{m})\right) \right\vert \right\vert ^{2}ds.
\end{eqnarray*}%
Applying Hölder's inequality and assumption $\left( H.11\right) $, we deduce
that%
\begin{eqnarray*}
&&\left( 1-\alpha \right) \left\{ \mathbb{E}\int_{t}^{T}\left\vert \tilde{Z}%
_{s}^{n}-\tilde{Z}_{s}^{m}\right\vert ^{2}ds+\mathbb{E}\int_{t}^{T}\int_{E}%
\left\vert \tilde{U}_{s}^{n}-\tilde{U}_{s}^{m}\right\vert ^{2}\lambda \left(
de\right)ds \right\}  \\
&\leq &2\mathbb{E}\left( \int_{t}^{T}\left\vert \tilde{Y}_{s}^{n}-\tilde{Y}%
_{s}^{m}\right\vert ^{2}ds\right) ^{\frac{1}{2}}\mathbb{E}\left(
\int_{t}^{T}\left\vert \Gamma _{s}^{n}-\Gamma _{s}^{m}\right\vert
^{2}ds\right) ^{\frac{1}{2}}+C\mathbb{E}\int_{t}^{T}\left\vert \tilde{Y}%
_{s}^{n}-\tilde{Y}_{s}^{m}\right\vert ^{2}ds.
\end{eqnarray*}%
The boundedness of the sequence $\left( \tilde{Y}^{n},\tilde{Z}^{n},\tilde{K}%
^{n},\tilde{U}^{n}\right) $, we deduce that $\Lambda =\sup_{n\in
\mathbb{N}
}\left( \mathbb{E}\int_{0}^{T}\left\vert \Gamma _{s}^{n}\right\vert
^{2}ds\right) <\infty ,$ this yields that%
\begin{eqnarray*}
&&\left( 1-\alpha \right) \mathbb{E}\int_{t}^{T}\left\vert \tilde{Z}_{s}^{n}-%
\tilde{Z}_{s}^{m}\right\vert ^{2}ds+\mathbb{E}\int_{t}^{T}\int_{E}\left\vert
\tilde{U}_{s}^{n}-\tilde{U}_{s}^{m}\right\vert ^{2}\lambda \left( de\right)
ds \\
&\leq &4\Lambda \mathbb{E}\left( \int_{t}^{T}\left\vert \tilde{Y}_{s}^{n}-%
\tilde{Y}_{s}^{m}\right\vert ^{2}ds\right) ^{\frac{1}{2}}+C\mathbb{E}%
\int_{t}^{T}\left\vert \tilde{Y}_{s}^{n}-\tilde{Y}_{s}^{m}\right\vert ^{2}ds.
\end{eqnarray*}%
Which yields that $\left( \tilde{Z}^{n}\right) _{n\geq 0}$ respectively $%
\left( \tilde{U}^{n}\right) _{n\geq 0}$ is a Cauchy sequence in $\mathcal{M}%
^{2}\left( 0,T,%
\mathbb{R}
^{d}\right) $ respectively in $\mathcal{L}^{2}\left( 0,T,\tilde{\mu},%
\mathbb{R}
\right) .$ Then there exists $\left( Z,U\right) \in \mathcal{M}^{2}\left(
0,T,%
\mathbb{R}
^{d}\right) \times \mathcal{L}^{2}\left( 0,T,\tilde{\mu},%
\mathbb{R}
\right) $ such that,%
\begin{equation}
\mathbb{E}\int_{t}^{T}\left\vert \tilde{Z}_{s}^{n}-Z_{s}\right\vert ^{2}ds+%
\mathbb{E}\int_{t}^{T}\int_{E}\left\vert \tilde{U}_{s}^{n}-U_{s}\right\vert
^{2}\lambda \left( de\right) \rightarrow 0,\text{\quad as }n\rightarrow
\infty .
\end{equation}%
On the other hand, applying Burkholder-Davis-Gundy inequality and $\left(
5.6\right) $, we obtain%
\begin{equation*}
\left\{
\begin{array}{l}
\mathbb{E}\sup_{0\leq t\leq T}\left\vert \int_{t}^{T}\tilde{Z}%
_{s}^{n}dW_{s}-\int_{t}^{T}Z_{s}dW_{s}\right\vert ^{2}\leq \mathbb{E}%
\int_{t}^{T}\left\vert \tilde{Z}_{s}^{n}-Z_{s}\right\vert ^{2}ds\rightarrow
0,\text{ as }n\rightarrow \infty , \\
\\
\mathbb{E}\sup_{0\leq t\leq T}\left\vert \int_{t}^{T}\int_{E}\tilde{U}%
_{s}^{n}\left( e\right) \tilde{\mu}\left( ds,de\right)
-\int_{t}^{T}\int_{E}U_{s}\left( e\right) \tilde{\mu}\left( ds,de\right)
\right\vert ^{2} \\
\leq \mathbb{E}\int_{t}^{T}\int_{E}\left\vert \tilde{U}_{s}^{n}-U_{s}\right%
\vert ^{2}\lambda \left( de\right) ds\rightarrow 0,\text{ as }n\rightarrow
\infty , \\
\\
\mathbb{E}\sup_{0\leq t\leq T}\left\vert \int_{t}^{T}g(s,\tilde{Y}_{s}^{n},%
\tilde{Z}_{s}^{n},\tilde{U}_{s}^{n})d\overleftarrow{B}_{s}-%
\int_{t}^{T}g(s,Y_{s},Z_{s},U_{s})d\overleftarrow{B}_{s}\right\vert ^{2} \\
\leq C\mathbb{E}\int_{t}^{T}\left\vert \tilde{Y}_{s}^{n}-Y_{s}\right\vert
^{2}ds+\alpha \mathbb{E}\int_{t}^{T}\left\vert \tilde{Z}_{s}^{n}-Z_{s}\right%
\vert ^{2}ds+\alpha \mathbb{E}\int_{t}^{T}\int_{E}\left\vert \tilde{U}%
_{s}^{n}-U_{s}\right\vert ^{2}\lambda \left( de\right) ds\rightarrow 0,\text{
as }n\rightarrow \infty .%
\end{array}%
\right.
\end{equation*}%
Therefore, from the properieties of $\left( f,\pi \right) $, we have%
\begin{equation*}
\Gamma _{s}^{n}=f(s,\tilde{Y}_{s}^{n-1},\tilde{Z}_{s}^{n-1},\tilde{U}%
_{s}^{n-1})+\pi \left( s,\delta \tilde{Y}_{s}^{n},\delta \tilde{Z}%
_{s}^{n},\delta \tilde{U}_{s}^{n}\right) \rightarrow f(s,Y_{s},Z_{s},U_{s}),
\end{equation*}%
$P-a.s.,$ for all $t\in \left[ 0,T\right] $ as $n\rightarrow \infty $. Then
follows by dominated convergence theorem that%
\begin{equation*}
\mathbb{E}\int_{0}^{T}\left\vert \Gamma
_{s}^{n}-f(s,Y_{s},Z_{s},U_{s})\right\vert ^{2}ds\rightarrow 0.
\end{equation*}%
Since $\left( \tilde{Y}_{s}^{n},\tilde{Z}_{s}^{n},\tilde{U}_{s}^{n},\Gamma
_{s}^{n}\right) $ converges in $\mathcal{B}^{2}\left(
\mathbb{R}
\right) \times \mathcal{M}^{2}\left( 0,T,%
\mathbb{R}
\right) $ and%
\begin{eqnarray*}
\mathbb{E}\left( \sup_{0\leq t\leq T}\left\vert \tilde{K}_{t}^{n}-\tilde{K}%
_{t}^{m}\right\vert ^{2}\right)  &\leq &\mathbb{E}\left\vert \tilde{Y}%
_{0}^{n}-\tilde{Y}_{0}^{m}\right\vert ^{2}+\mathbb{E}\left( \sup_{0\leq
t\leq T}\left\vert \tilde{Y}_{t}^{n}-\tilde{Y}_{t}^{m}\right\vert
^{2}\right) +\mathbb{E}\int_{0}^{T}\left\vert \Gamma _{s}^{n}-\Gamma
_{s}^{m}\right\vert ^{2}ds \\
&&+\mathbb{E}\sup_{0\leq t\leq T}\left\vert \int_{0}^{t}\left( g(s,\tilde{Y}%
_{s}^{n},\tilde{Z}_{s}^{n},\tilde{U}_{s}^{n})-g(s,\tilde{Y}_{s}^{m},\tilde{Z}%
_{s}^{m},\tilde{U}_{s}^{m})\right) d\overleftarrow{B}\right\vert ^{2} \\
&&+\mathbb{E}\sup_{0\leq t\leq T}\left\vert \int_{0}^{t}\left( \tilde{Z}%
_{s}^{n}-\tilde{Z}_{s}^{m}\right) dW_{s}\right\vert ^{2}+\mathbb{E}%
\sup_{0\leq t\leq T}\left\vert \int_{0}^{t}\int_{E}\left( \tilde{U}%
_{s}^{n}\left( e\right) -\tilde{U}_{s}^{m}\left( e\right) \right) \tilde{\mu}%
\left( ds,de\right) \right\vert ^{2},
\end{eqnarray*}%
for any $n,m\geq 0$, we deduce that%
\begin{equation*}
\mathbb{E}\left( \sup_{0\leq t\leq T}\left\vert \tilde{K}_{t}^{n}-\tilde{K}%
_{t}^{m}\right\vert ^{2}\right) \rightarrow 0,
\end{equation*}%
as $n,m\rightarrow \infty .$ Consequently, there exists a $\mathcal{F}_{t}-$%
measurable process $K$ wich value in $%
\mathbb{R}
$ such that%
\begin{equation}
\mathbb{E}\left( \sup_{0\leq t\leq T}\left\vert \tilde{K}_{t}^{n}-K_{t}%
\right\vert ^{2}\right) \rightarrow 0,\text{\quad }as\text{ }n\rightarrow
\infty .
\end{equation}%
Finally, we have%
\begin{equation*}
\mathbb{E}\left( \sup_{0\leq t\leq T}\left\vert \tilde{Y}_{t}^{n}-Y_{t}%
\right\vert ^{2}+\int_{t}^{T}\left\vert \tilde{Z}_{s}^{n}-Z_{s}\right\vert
^{2}ds+\int_{t}^{T}\int_{E}\left\vert \tilde{U}_{s}^{n}-U_{s}\right\vert
^{2}\lambda \left( de\right)ds +\sup_{0\leq t\leq T}\left\vert \tilde{K}%
_{t}^{n}-K_{t}\right\vert ^{2}\right) \rightarrow 0,\text{\quad }as\text{ }%
n\rightarrow \infty .
\end{equation*}%
Obviously, $K_{0}=0$ and $\{K_{t};0\leq t\leq T\}$ is a increasing
and
continuous process. From $(5.4)$, we have for all $n\geq 0$, $\tilde{Y}%
_{t}^{n}\geq S_{t}$, $\forall t\in \lbrack 0,T]$, then $Y_{t}\geq S_{t}$, $%
\forall t\in \lbrack 0,T]$.

On the other hand, from the result of Saisho, we have%
\begin{equation*}
\int_{0}^{T}\left( \tilde{Y}_{s}^{n}-S_{s}\right) d\tilde{K}%
_{s}^{n}\rightarrow \int_{0}^{T}\left( Y_{s}-S_{s}\right) dK_{s},\text{ }%
\mathbf{P}-a.s.,\text{\quad }as\text{ }n\rightarrow \infty .
\end{equation*}%
Using the identity $\int_{0}^{T}\left( \tilde{Y}_{s}^{n}-S_{s}\right) d%
\tilde{K}_{s}^{n}=0$ for all $n\geq 0,$ we obtain
$\int_{0}^{T}\left( Y_{s}-S_{s}\right) dK_{s}=0.$ letting
$n\rightarrow +\infty $ in equation $\left( 1.1\right) $, we prove
that $\left( Y_{t},Z_{t},K_{t},U_{t}\right) _{t\in \lbrack 0,T]}$ is
the solution to $\left( 1.1\right) .$

Let $\left( Y_{\ast },Z_{\ast },U_{\ast },K_{\ast }\right) $ be a
solution of $\left( 1.1\right) $. Then by Theorem 3.1, we have for
any $n\in
\mathbb{N}
^{\ast }$, $Y^{n}\leq Y_{\ast }.$ Therefore, $Y$ is a minimal solution of $%
\left( 1.1\right) .$ \eop
\begin{remark}
Using the same arguments and the following approximating sequence
\begin{equation*}
f_{n}\left( t,\omega ,y,z,u\right) =\sup_{\left( y^{^{\prime
}},z^{\prime },u^{^{\prime }}\right) \in \mathcal{B}^{2}\left(
\mathbb{R}
\right) }\left[ f\left( t,\omega ,y^{^{\prime }},z^{\prime
},u^{^{\prime }}\right) -n\left( \left\vert y-y^{^{\prime
}}\right\vert +\left\vert
z-z^{\prime }\right\vert +\left\vert u-u^{^{\prime }}\right\vert \right) %
\right] ,
\end{equation*}%
one can prove that the RBDSDE (1.1) has a maximal solution.
\end{remark}

\end{document}